\documentclass[a4paper]{amsart}
\date{Jume 02, 2003}
\usepackage[dvips]{graphicx}
\usepackage{verbatim,enumerate}
\usepackage{amsthm} 
\theoremstyle{plain}
 \newtheorem{theorem}{Theorem}[section]
 \newtheorem*{theorem*}{Theorem}
 \newtheorem*{proposition*}{Proposition}
 \newtheorem{proposition}[theorem]{Proposition}
 \newtheorem{lemma}[theorem]{Lemma}
 \newtheorem{corollary}[theorem]{Corollary}
\theoremstyle{remark}
 \newtheorem{definition}[theorem]{Definition}
 \newtheorem{remark}[theorem]{Remark}
 \newtheorem*{remark*}{Remark}
 \newtheorem{example}[theorem]{Example}
\numberwithin{equation}{section}

\newcommand{\Z}{\boldsymbol{Z}}
\newcommand{\R}{\boldsymbol{R}}
\newcommand{\C}{\boldsymbol{C}}
\newcommand{\CP}{\boldsymbol{C\!P}}

\newcommand{\trace}{\operatorname{trace}}
\newcommand{\ord}{\operatorname{ord}}

\newcommand{\SL}{\operatorname{SL}}
\newcommand{\SU}{\operatorname{SU}}
\newcommand{\U}{\operatorname{U}}
\newcommand{\PSU}{\operatorname{PSU}}
\newcommand{\PSL}{\operatorname{PSL}}
\newcommand{\SO}{\operatorname{SO}}
\newcommand{\Herm}{\operatorname{Herm}}
\newcommand{\F}{\mathcal{F}}
\newcommand{\G}{\mathcal{G}}
\newcommand{\zb}{\bar z}
\renewcommand{\Re}{\operatorname{Re}}
\newcommand{\Geod}{\operatorname{Geod}}
\title[Flat Fronts]{
   Flat fronts in hyperbolic 3-space
}
\author{Masatoshi Kokubu}
\address[Kokubu]{%
   Department of Natural Science,
   School of Engineering,
   Tokyo Denki University,
   2-2 Kanda-Nishiki-Cho,
   Chiyoda-Ku, Tokyo, 101-8457,
   Japan
}
\email{kokubu@cck.dendai.ac.jp}
\author{Masaaki Umehara}
\address[Umehara]{%
   Department of Mathematics, Graduate School of Science,
   Osaka University,
   Toyonaka, Osaka 560-0043,
   Japan
}
\email{umehara@math.wani.osaka-u.ac.jp}
\author{Kotaro Yamada}
\address[Yamada]{%
   Faculty of Mathematics,
   Kyushu University 36, 
   Higashi-ku, Fukuoka 812-8581, Japan%
}
\email{kotaro@math.kyushu-u.ac.jp}
\subjclass{Primary 53C42; Secondary 53A35, 53D99}
\begin{document}
\begin{abstract}
 We shall investigate flat surfaces in hyperbolic 3-space 
 with admissible singularities, called `flat fronts'. 
 An Osserman-type inequality for complete flat fronts is shown. 
 When equality holds in this inequality, we show that all 
 the ends are embedded, 
 and give new examples for which equality holds.
\end{abstract}
\maketitle
\section*{Introduction}
It is a classical fact that any complete flat surface in the hyperbolic
$3$-space $H^3$ must be a horosphere or a hyperbolic cylinder.
However, this does not imply the lack of an interesting global theory
for flat surfaces.
Recently, G\'alvez, Mart\'\i{}nez and Mil\'an \cite{GMM} established a
Weierstrass-type representation formula for such surfaces. 
More recently, the authors \cite{KUY} proved another representation
formula constructing flat surface from a given pair of hyperbolic Gauss
maps, and also gave new examples.

In this paper, we shall investigate global properties of flat surfaces
with `admissible'  singularities, 
which contains all of the previous examples in \cite{GMM} and
\cite{KUY}.
(A singular (i.e.\ degenerate) point is called admissible if the
corresponding points on nearby parallel surfaces are regularly
immersed. See Section~\ref{sec:front}.)
Such surfaces are characterized as the projections of Legendrian
immersions in the unit cotangent bundle $T_1^* H^3$ of $H^3$, called
{\it flat fronts}:
The $5$-manifold $T_1^* H^3$ has the canonical contact form $\eta$.
If we identify $H^3$ with the Poincar\'e ball $(D^3;x^1,x^2,x^3)$, any
element $\alpha$ of the cotangent bundle $T^*H^3$ can be written as
\[
   \alpha=p_1(\alpha)\, dx^1 + p_2(\alpha)\, dx^2 + 
                 p_3(\alpha)\, dx^3 \qquad(\in T^*H^3).
\]
Then $(p_1,p_2,p_3,x^1,x^2,x^3)$ gives a canonical coordinate system of
$T^*H^3$ and the canonical form on $T^*H^3$ 
\[
       \eta=p_1\, dx^1 + p_2\, dx^2 + p_3\, dx^3, 
\]
which induces a canonical contact form on $T^*_1 H^3$.
A map $L \colon M^2\to T_1^* H^3$ is called {\it Legendrian\/} if
the pull-back $L^*\eta$  vanishes identically.
For a given immersion $f \colon M^2\to H^3$, there exists a unique
Legendrian immersion
\[
    L_f \colon M^2\longrightarrow T_1^* H^3
\]
such that $\pi\circ L_f=f$, 
where $\pi \colon T_1^* H^3\to H^3$ is the projection. 
That is, any immersion can be lifted to a Legendrian immersion. 
However, the converse is not true.
A projection 
\[
         \pi\circ L \colon M^2 \longrightarrow H^3
\]
of a Legendrian immersion $L$ is called a {\it {\rm (}wave{\rm )} front\/}, 
which may have singular points 
(i.e. the point where the Jacobi matrix degenerates.) 
A point which is not singular is called {\it regular\/}, 
where the first fundamental form is positive definite.
The Gaussian curvature is well-defined at regular points.
A front is called {\it flat\/} if the Gaussian curvature vanishes at
each regular point. 

A front $f$ is called {\it complete\/} if there is a symmetric tensor
$T$ on $M^2$ which has compact support such that $T+ds^2$ is a complete
Riemannian metric on $M^2$,
where $ds^2$ is the first fundamental form of $f$.
If  $M^2$ is orientable,
$M^2$ can be regarded as a Riemann surface whose complex structure
is compatible with  respect to the pull-back of the Sasakian metric on
$T_1^*H^3$ by $L_f$.
Moreover, the second fundamental form is hermitian with respect to this
structure, and there is a closed Riemann surface $\overline M^2$
such that $M^2$ is biholomorphic to 
$\overline M^2\setminus \{p_1,\dots,p_n\}$.
The points $p_1,\dots,p_n$ are called the {\it ends\/} of $f$.

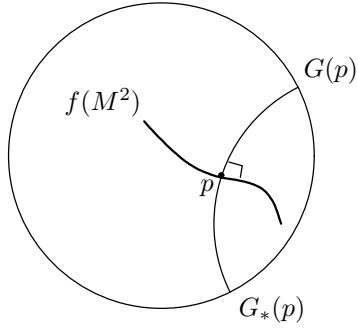
\begin{figure}
\unitlength 0.1in
\begin{picture}( 16.0000, 16.0000)(  7.5800,-19.2700)
%
\special{pn 8}%
\special{ar 1558 1128 800 800  0.0000000 6.2831853}%
%
\special{pn 8}%
\special{ar 2632 1484 800 800  2.6773858 4.2481821}%
%
\special{pn 13}%
\special{pa 2184 1484}%
\special{pa 2174 1452}%
\special{pa 2162 1420}%
\special{pa 2150 1388}%
\special{pa 2136 1360}%
\special{pa 2118 1334}%
\special{pa 2100 1312}%
\special{pa 2076 1292}%
\special{pa 2050 1278}%
\special{pa 2022 1266}%
\special{pa 1992 1256}%
\special{pa 1960 1248}%
\special{pa 1928 1240}%
\special{pa 1894 1234}%
\special{pa 1860 1226}%
\special{pa 1828 1218}%
\special{pa 1798 1206}%
\special{pa 1768 1194}%
\special{pa 1738 1180}%
\special{pa 1710 1164}%
\special{pa 1684 1146}%
\special{pa 1658 1128}%
\special{pa 1632 1108}%
\special{pa 1608 1088}%
\special{pa 1584 1066}%
\special{pa 1560 1044}%
\special{pa 1538 1020}%
\special{pa 1514 998}%
\special{pa 1492 974}%
\special{pa 1470 950}%
\special{pa 1468 948}%
\special{sp}%
\put(14.5000,-9.2000){\makebox(0,0)[rb]{$f(M^2)$}}%
%
\special{pn 8}%
\special{pa 1910 1160}%
\special{pa 1980 1180}%
\special{fp}%
\special{pa 1980 1180}%
\special{pa 1970 1240}%
\special{fp}%
%
\special{pn 13}%
\special{sh 1}%
\special{ar 1870 1230 10 10 0  6.28318530717959E+0000}%
\special{sh 1}%
\special{ar 1870 1230 10 10 0  6.28318530717959E+0000}%
\put(18.3000,-12.5000){\makebox(0,0)[rt]{$p$}}%
\put(23.0000,-7.5000){\makebox(0,0)[lb]{$G(p)$}}%
\put(19.6000,-18.6000){\makebox(0,0)[lt]{$G_{*}(p)$}}%
\end{picture}%
\caption{Hyperbolic Gauss maps}
\label{fig:hgmap}
\end{figure}

For each point $p\in M^2$, there exists a pair $(G(p),G_*(p))
\in S^2\times S^2$ of distinct points in the ideal boundary
$S^2=\partial H^3$  such that the geodesic in $H^3$ starting from
$G_*(p)$ towards  $G(p)$ coincides the oriented
normal geodesic at $p$ (see Figure~\ref{fig:hgmap}).
The maps 
\[
   G, ~ G_*\colon \overline M^2\setminus \{p_1,\dots,p_n\} 
      \longrightarrow S^2 
\]
are called the {\it positive\/} and 
{\it negative hyperbolic Gauss maps\/} of $f$, respectively.
They are holomorphic if we regard $S^2=\partial H^3$ as the Riemann
sphere.
An end $p_j$ is called {\it regular\/} if both of $G$ and $G_*$
extend holomorphically across it.
As we shall show later, there are many flat fronts with regular ends.
Moreover, such surfaces
satisfy the following global property:

\begin{theorem*}
 An orientable
 complete flat front 
 $f\colon \overline M^2\setminus \{p_1,\dots,p_n\}\to H^3$ 
 with regular ends satisfies the inequality
 \[
       \deg G+\deg G_*\ge n.
 \]
 Moreover, equality holds if and only if all ends are embedded.
\end{theorem*}
Here, $\deg G$ denotes the degree of a holomorphic map 
$G\colon{}\overline M^2\to\CP^1=S^2$.
This inequality is an analogue of (resp.\ a hyperbolic version of) the
Osserman inequality
\[
       2\deg G+\chi\bigl(\overline M^2\setminus \{p_1,\dots,p_n\}\bigr) 
       \ge n
\]
for a complete minimal (resp.\ a mean curvature one) surface 
$f\colon \overline M^2\setminus \linebreak\{p_1,\dots,p_n\}\to \R^3$
(resp.\ $H^3$) with finite total curvature, 
where $G$ is the Gauss map (resp.\ hyperbolic Gauss map) of the surface.
In these two cases, like as ours, equality implies 
the embeddedness of ends
(see \cite{Oss,JM} for minimal surface case, and  \cite{UY} for mean
curvature one case).

To prove that equality implies the ends are embedded,
a criterion for embeddedness of ends given in \cite{GMM} will be applied.
Furthermore, we shall classify flat $3$-noids and give 
a genus $1$ flat front with regular ends
(Section~\ref{sec:example}). 

\bigskip
On the other hand, in Section~\ref{sec:complete}, we discuss
an alternative global inequality for flat fronts:
Let $d\sigma^2$ be the pull-back of the Sasakian metric by
the Legendrian lift of a complete flat front $f$
(see \eqref{eq:sasakian} in Section~\ref{sec:front}).
Then $d\sigma^2$ is complete and the Gaussian curvature
$K_{d\sigma^2}$ is non-positive.
Then we show a geometric condition for equality of the Cohn-Vossen
inequality:
\begin{proposition*}
 For a complete flat front $f\colon{}M^2\to H^3$, it holds that
 \[
    \frac{1}{2\pi}\int_{M^2}(-K_{d\sigma^2})\,dA_{d\sigma^2}
      \geq -\chi(M^2),
 \]
 where $dA_{d\sigma^2}$ is the area element of $d\sigma^2$ and
 $\chi(M^2)$ is the Euler number of $M^2$.
 Moreover, the equality holds if and only if all ends are asymptotic
 to a hyperbolic cylinder.
\end{proposition*}
Note that flat hypersurfaces in $H^n$ ($n\geq 4$)  are totally
umbilic.
So $n=3$ is the interesting case.
\bigskip

The authors are very grateful to P. Pirola and E. Musso for their
valuable comments.

\section{Local properties of flat surfaces}
\label{sec:preliminaries}
In this section, we review local properties of flat immersions.
We denote by $H^3$ the hyperbolic 3-space of constant curvature $-1$.
Let $M^2$ be a 2-manifold and 
\[
   f \colon M^2 \longrightarrow H^3
\]
be a flat immersion, that is, the Gaussian curvature of the induced
metric vanishes.
Then it follows from the Gauss equation that the second fundamental 
form is positive or negative definite and thus $M^2$ is 
orientable.
We fix an orientation of $M^2$, then $M^2$ can be regarded
as a Riemann surface such that the second fundamental
form $dh^2$ is hermitian.
A holomorphic (resp.\ meromorphic) map 
\[
   E=
     \begin{pmatrix}
      A & B \\ C & D
     \end{pmatrix}%
     \colon M^2\longrightarrow \PSL(2,\C)
\]
is called {\it Legendrian\/} if
\begin{equation}\label{cnd:legendre}
       D\,dA-B\, dC=0
\end{equation}
holds. Indeed, \eqref{cnd:legendre} implies that the vanishing of 
the pull-back of a holomorphic contact form on $\PSL(2, \C)$.   
As is shown in \cite{GMM}, there exists a holomorphic Legendrian
immersion $E_f$ (called a {\it holomorphic Legendrian lift\/} of $f$)
defined on the universal cover $\widetilde M^2$ of $M^2$
\[
    E_f\colon \widetilde M^2 \longrightarrow \PSL(2,\C)
\]
such that $f$ is the projection of $E_f$ onto $H^3=\PSL(2,\C)/\PSU(2)$.
Since \eqref{cnd:legendre} implies $E_f^{-1}dE_f$ is off-diagonal, 
we can set 
\[
     E_f^{-1}dE_f
       =\begin{pmatrix}
	 0 & \theta \\
	 \omega & 0
      \end{pmatrix}.
\]
The holomorphic $1$-forms $\omega$ and $\theta$ are called the 
{\it first canonical form\/} and the {\it second canonical form\/}, 
respectively.
It holds that
\begin{align}
  \label{eq:omega-def}
    \omega&=\begin{cases}
	     \dfrac{dA}{B}
	     \qquad
	     &
	     \bigl(\text{if $dA\not\equiv0$ or $B\not\equiv 0$}\bigr),\\[8pt]
	     \dfrac{dC}{D}
	     &\bigl(\text{if $dC\not\equiv0$ or $D\not\equiv 0$}\bigr),\\[8pt]
	  \end{cases}\\
 \label{eq:theta-def}
    \theta&=\begin{cases}
	     \dfrac{dB}{A}
	     \qquad
             & 
	     \bigl(\text{if $dB\not\equiv0$ or $A\not\equiv 0$}\bigr),\\[8pt]
	     \dfrac{dD}{C}
	     &\bigl(\text{if $dD\not\equiv0$ or $C\not\equiv 0$}\bigr).\\[8pt]
	  \end{cases}
 \end{align}
Here $dA\not\equiv 0$ (resp.~$B\not\equiv 0$) means that the $1$-form
$dA$ (resp.~the function $B$) is not identically zero.

In particular,  if all cases in \eqref{eq:omega-def} and
\eqref{eq:theta-def} are well-defined,
\begin{equation}\label{eq:omega-theta}
     \omega= \frac{dA}{B}=\frac{dC}{D} \qquad\text{and}\qquad
     \theta= \frac{dB}{A}=\frac{dD}{C}
\end{equation}
holds.
Then the first and second fundamental forms $ds^2$ and $dh^2$ have the
following expressions
\begin{align}
  ds^2&=(\omega+\bar \theta)(\bar \omega+ \theta) 
          =\omega\theta +\bar \omega \bar \theta+
                    |\omega|^2+|\theta|^2 ,\label{eq:first}\\
  dh^2&=|\theta|^2 - |\omega|^2.\label{eq:second}
\end{align}
 Though $\omega$ and $\theta$ are defined only on the universal 
 cover $\widetilde M^2$, the first fundamental form $ds^2$
 is well-defined on $M^2$, and then so is the $(1,1)$-part of $ds^2$:
\begin{equation}\label{eq:one-one-part}
   ds^2_{1,1}:=|\omega|^2+|\theta|^2.
\end{equation}
Since \eqref{eq:second} is well-defined on $M^2$, 
so are $|\omega|^2$ and $|\theta|^2$.
Moreover, we can deduce that
\begin{equation}\label{eq:flat-pseudo}
\text{
\begin{tabular}{p{10cm}}
 $|\omega|^2$ and $|\theta|^2$ define flat pseudometrics on $M^2$
 which are compatible with the complex structure of $M^2$.
\end{tabular}}
\end{equation}

The $(2,0)$-part of $ds^2$ is called the 
{\it Hopf differential\/}, and it is denoted by $Q$, that is, 
\begin{equation}\label{def_hopf}
   Q:=\omega \theta.
\end{equation}

The {\it positive hyperbolic Gauss map\/} $G$ and the 
{\it negative hyperbolic Gauss map\/} $G_*$ of the flat surface
are defined by 
\begin{equation}\label{def_Gauss}
   G=\frac{A}{C},\qquad
   G_*=\frac{B}{D}. 
\end{equation}
They are single-valued on $M^2$. 
The geometric meaning of $G$ and $G_*$ are described in Introduction. 
(See also \cite{GMM}.)
By definition, 
\begin{equation}\label{eq:omega}
   dG= d \left( \frac{A}{C}\right)
     =\frac{dA\,C-A\, dC}{C^2}
     =\frac{BC-DA}{C^2}\omega=-\frac{\omega}{C^2}.
\end{equation}
Similarly, we have
\begin{equation}\label{eq:theta}
  dG_*=d \left( \frac{B}{D}\right)
      =\frac{dB\,D-B\, dD}{D^2}
      =\frac{AD-BC}{D^2}\theta=\frac{\theta}{D^2}.
\end{equation}
On the other hand, we have
\begin{equation}\label{eq:difference}
    G-G_*=\frac{A}{C}-\frac{B}{D}
         =\frac{AD-BC}{CD}=\frac{1}{CD}.
\end{equation}
We have the following identity
\begin{equation}\label{eq:Hopf}
     Q=\omega\theta=-(CD)^2dGdG_*=-\frac{dGdG_*}{(G-G_*)^2}.
\end{equation}

Now we set
\begin{equation}\label{eq:small-g}
       g(q):=\int_{p_0}^q \omega, \qquad
       g_*(q):=\int_{p_0}^q \theta
         \qquad (q\in M^2)
\end{equation}
where $p_0$ is a base point.
Then $g$ and $g_*$ are holomorphic functions defined on $\widetilde M^2$.
We remark that $(g,G)$ and $(g_*,G_*)$ satisfy the following
important relation (see \cite{GMM}):
\begin{equation}\label{eq_sch}
         S(g)-S(G)=2Q,\qquad S(g_*)-S(G_*)=2Q,
\end{equation}
where  $S(G)$ is the Schwarzian derivative 
\[
  S(G)=
  \left[
  \left(\frac{G''}{G'}\right)'-\frac 12 \left(\frac{G''}{G'}\right)^2
  \right]\,dz^2 \qquad
  \left( '=\frac{d}{dz}\right)
\]
with respect to a local complex coordinate $z$ on $M^2$.
Though the meromorphic $2$-differentials $S(g)$ and $S(G)$ depend on
complex coordinates, the difference $S(g)-S(G)$ does not.

\begin{remark}
 Hyperbolic $3$-space $H^3$ can be realized as a
 hyperboloid in Minkowski 4-space $(L^4,(x^0,x^1,x^2,x^3))$:
 \[
   H^3=\left\{
          (x^0,x^1,x^2,x^3)\in L^4\,;\, x^0>0,\,\, 
                   -(x^0)^2+\sum_{j=1}^3 (x^j)^2=-1
       \right\}.
 \]
 Let $f\colon M^2\to H^3$ be a flat immersion and assume $M^2$ is
 connected.
 Then the universal cover $\widetilde M^2$ of $M^2$ is diffeomorphic to
 $\R^2$ and has a coordinate system $(x,y)$ defined on $\widetilde M^2$
 such that the first fundamental form $ds^2$ can be written as
 \[
    ds^2=dx^2 + dy^2.
 \]
 Then we have an orthonormal frame field 
 \[
    e\colon \widetilde M^2 \ni p\mapsto 
     \bigl(f(p),f_x(p),f_y(p),\nu(p)\bigr)\in \SO(3,1),
 \]
 where $\nu(p)\in T_pH^3(\subset L^4)$ is a unit normal vector 
 of the immersion $f$ at $p$.
 Now, we can identify $L^4$ with the set $\Herm(2)$  of 2 by 2 Hermitian
 matrices, that is, 
 \begin{equation}\label{eq:lor-herm}
   L^4\ni (x^0,x^1,x^2,x^3) \longleftrightarrow
   \begin{pmatrix}
      x^0+x^3 & x^1+i x^2 \\
      x^1-ix^2 & x^0-x^3
   \end{pmatrix}
   \in \Herm(2).
 \end{equation}
 Then the hyperbolic 3-space $H^3$ can be rewritten as
 \begin{align*}
   H^3&=\{X\in \Herm(2)\,;\, \det(X)=1, \,\,\trace X>0\}\\
      &=\{aa^*\,;\, a\in\SL(2,\C)\},
 \end{align*}
 where $a^*={}^t\bar a$.
 Setting
 \[
      v_0:=\begin{pmatrix}
	    1 & 0 \\
	    0 & 1
	   \end{pmatrix}, \quad
      v_1:=\begin{pmatrix}
	    0 & i \\
	    i & 0
	   \end{pmatrix}, \quad
      v_2:=\begin{pmatrix}
	    \hphantom{-}0 & i \\
	    -i & 0
	   \end{pmatrix}, \quad
      v_3:=\begin{pmatrix}
	    1 & \hphantom{-}0 \\
	    0 & -1
	   \end{pmatrix}, 
 \]
 then there is a lift $E\colon \widetilde M^2\to \SL(2,\C)$ 
 of the orthonormal frame $e$ such that $e=\pi\circ E$ where 
 $\pi\colon \SL(2,\C)\to \SO(3,1)$ is the $2$-fold covering homomorphism,
 that is,
 \begin{equation}
  f=EE^*,\quad f_x=Ev_1E^*,\quad
   f_y=Ev_2E^*,\quad \nu=Ev_3E^*.
 \end{equation}
 Thus $E$ coincides with $E_f$. 
 This implies that $E$ itself is holomorphic with respect
 to the complex structure induced from the second fundamental form.
 A matrix multiplication $E\to aE$ ($a\in\SL(2,\C)$) corresponds
 to an  isometric change of the surface $f\mapsto a f a^*$.
 This induces the change of hyperbolic Gauss maps as
 \begin{equation}\label{eq:gauss-change}
    G \mapsto a\star G := \frac{a_{11}G+a_{12}}{a_{21}G+a_{22}},\qquad
    G_* \mapsto a\star G_* = \frac{a_{11}G_*+a_{12}}{a_{21}G_*+a_{22}},
 \end{equation}
 where $a=(a_{ij})$.

 It is interesting to compare this with the case of constant mean
 curvature one surfaces in $H^3$. 
 In that case, there is a holomorphic immersion 
 $F\colon \widetilde M^2\to \SL(2,\C)$ such that $f=FF^*$, 
 but it does not coincide with the lift 
 $E\colon \widetilde M^2 \to \SL(2,\C)$ 
 of an orthonormal frame.
 We must adjust $E$ by multiplying by a local $\SU(2)$-section 
 $s\colon\widetilde M^2\to \SU(2)$ so that $F:=Es$ becomes
 holomorphic (see Bryant~\cite{Br}).
\end{remark}
\section{Flat surfaces as (wave) fronts}
\label{sec:front}

In this section, we shall define `flat fronts' as projections of
Legendrian immersions into the unit cotangent bundle 
$T_1^*H^3$. 
Since $T_1^*H^3$ is isomorphic to the unit tangent bundle $T_1H^3$, 
we can make the following identification
\[
   T_1^*H^3\cong \F
         :=\bigl\{(x,v)\in L^4 \times L^4\,;\,
	           -\langle x,x\rangle=
         	   \langle v, v\rangle=1,\langle x,v\rangle=0
	   \},
\]
where $\langle~,~\rangle$ is the inner product of $L^4$.
The metric 
\[
    d\sigma_0^2:=\sum_{j=0}^3 (dx^j)^2
           +\sum_{j=0}^3 (dv^j)^2
	   \quad 
	   \bigl(x=(x^0,x^1,x^2,x^3),~
	         v=(v^0,v^1,v^2,v^3)\bigr)
\]
on $\F$ induced from the product of Lorentzian metrics of $L^4\times
L^4$ is positive definite, and is called the {\it Sasakian metric}.
In fact, if we identify $\F$ with $T_1H^3$,  it coincides with the metric
on the unit tangent bundle defined by Sasaki \cite{S1,S2}.
The contact form of $\F$ is given by
\[
       \eta:=\sum_{j=0}^3 v^jdx^j.
\]
Now, a Legendrian map $L$ of a $2$-manifold $M^2$ into the unit
cotangent bundle can be identified with a map
\[
     L\colon M^2\longrightarrow \F
\]
such that $L^*\eta$ vanishes. 
We denote two canonical projections by
\[
    \pi_F\colon \mathcal F\ni (x,v)\longmapsto x \in H^3,\qquad
    \pi'_F\colon \mathcal F\ni (x,v)\longmapsto v \in L^4.
\]
A map $f\colon M^2\to H^3$ is called a {\it front\/} if there exists
a Legendrian immersion $L_f\colon M^2\to \F$ called the 
{\it Legendrian lift\/} of $f$ such that
\[
     \pi_F\circ L_f =f.
\]

We remark that, by definition, any immersion $f \colon M^2\to H^3$ is a
front if $M^2$ is orientable. 
In fact, $L_f$ is given by a pair $(f,\nu_f)$ of $f$ and the unit normal
vector $\nu_f$ of $f$.

For a given front $f\colon M^2\to H^3$, 
we can define a {\it parallel front\/} $f_t\colon M^2\to H^3$ of
distance $t$ as follows
\[
     f_t:=(\cosh t)f+ (\sinh t) \nu_f=\pi_F\circ L_t,
\]
where 
\[
   L_t:=(f_t,\nu^{}_{f_t})\qquad
   \bigl(\nu^{}_{f_t}:=(\sinh t)f + (\cosh t)\nu_f\bigr)
\]
is a Legendrian immersion and
\[
    \nu_f:=\pi'_F\circ L_f\colon M^2\longrightarrow L^4.
\]

When $f$ is an immersion, this is nothing but the definition of a
parallel surface. 
So we call $\nu_f$ the unit normal vector (field) of the front $f$.
 
For a given front $f\colon M^2\to H^3$, 
\[
   ds^2:=\langle df, df \rangle \qquad \text{and}\qquad
   dh^2:=-\langle df, d\nu_f \rangle
\]
are called {\em the first\/} and {\em the second fundamental forms},
respectively.

\begin{definition}
 A front $f\colon{}M^2\to H^3$ is called {\it flat\/} 
 if for each $p\in M^2$, there exists a real number $t\in \R$
 such that the parallel front $f_t$ gives a flat immersion at $p$. 
\end{definition}

\begin{remark}
 An equivalent definition of a flat front is that the Gaussian
 curvature of $f$ at the regular points vanishes. 
 However, this definition is not suitable when all points of
 $f$ degenerate, 
 and such a case really occurs,
 since hyperbolic cylinders can collapse to a geodesic.
\end{remark}

As shown in the following proposition,
all parallel fronts $f_t$ ($t\in \R$) of a flat 
front $f$ are also flat fronts.

\begin{proposition}
\label{prop:parallel}
 Let $f\colon M^2\to H^3$ be a flat front, then the second fundamental
 form $dh^2$ is proportional to the pull-back of the Sasakian  metric
 $d\sigma^2=L_f^*d\sigma^2_0$. 
 Moreover, the parallel front $f_t$ of $f$ is also a flat front for all
 $t$. 
 In particular, the Gaussian curvature of $f$ at the regular point
 vanishes.
\end{proposition}

\begin{remark}
As in \cite{GMM}, the lift $E_{f_t}$ of $f_t$ is 
given by
\[
  E_{f_t}=E_{f}
    \begin{pmatrix}
     e^{-t/2} & 0 \\
     0 & e^{t/2}
    \end{pmatrix}.
\]
\end{remark}
\begin{proof}[Proof of Proposition~\ref{prop:parallel}]
We fix a point $p\in M^2$.
Then, by definition, there is a parallel front 
$f_{t_0}\colon M^2\to H^3$ such that $f_{t_0}$ is regular at $p$ 
and the Gaussian curvature of $f_{t_0}$ vanishes around $p$.
Without loss of generality, 
we may assume that $t_0=0$, that is, $f=f_{t_0}$.

First, we consider the case that 
the first and second fundamental forms are proportional.
Then $f$ must be a horosphere and the statement of the theorem is
obvious.

So we may assume that the second fundamental form is not proportional to
the first fundamental form.
We can write the Legendrian lift $L_f$ as a pair $L_f=(f,\nu_f)$,
where $\nu_f$ is the unit normal vector field of $f$.
Then we have
\begin{equation}\label{eq:sasakian}
   d\sigma^2=
        \langle df, df \rangle + \langle d\nu, d\nu \rangle.
\end{equation}
Now we fix a local coordinate neighborhood
$(U;u,v)$ of $M^2$ and define three 2 by 2 matrices
as follows:
\begin{align*}
  M_1
    &:=\hphantom{-}\begin{pmatrix}
	\langle f_u, f_u \rangle & \langle f_u, f_v \rangle \\
	\langle f_u, f_v \rangle & \langle f_v, f_v \rangle \\
       \end{pmatrix}, \\
  M_2
    &:=-\begin{pmatrix}
	 \langle f_u, \nu_u \rangle & \langle f_u, \nu_v \rangle \\
	 \langle f_u, \nu_v \rangle & \langle f_v, \nu_v \rangle \\
	\end{pmatrix}, \\
  M_3
    &:=\hphantom{-}\begin{pmatrix}
	\langle \nu_u, \nu_u \rangle & \langle \nu_u, \nu_v \rangle \\
	\langle \nu_u, \nu_v \rangle & \langle \nu_v, \nu_v \rangle \\
       \end{pmatrix} .
\end{align*}
We set $A_f:=M_1^{-1}M_2$, which is the shape operator of $f$.
The Gauss equation implies that
\begin{equation}\label{eq:Gauss}
   \det A_f =1+K_{ds^2},
\end{equation}
where $K_{ds^2}$ is the Gaussian curvature of $ds^2$.
On the other hand, by the definition of 
$M_3$, we have
\[
       M_3=M_2A=M_1A^2. 
\]
By the Cayley-Hamilton theorem, we have
\[
   M_3=M_1\bigl(2(\trace A)A-(1+K_{ds^2})I\bigr)=
         2(\trace A)M_2-K_{ds^2}M_1,
\] 
where $I$ is the identity matrix. 
Thus we have
\[
   M_1+M_3-2(\trace A) M_2=-K_{ds^2}M_1.
\]
Since $M_1$ is not proportional to $M_2$,
this implies that $M_1+M_3$ is proportional to $M_2$ if and only if
$K_{ds^2}$ vanishes.
So the second fundamental form $dh^2$ is proportional to  $d\sigma^2$
when $f$ is flat.
Now we shall show that $f_t$ is also flat.
In fact, $f_t$ and its unit normal vector $\nu_t$ have the following
expressions
\[
       f_t=(\cosh t)f+(\sinh t)\nu,
             \qquad
       \nu_t=(\sinh t)f+(\cosh t)\nu.
\]
The fundamental forms are given by
\begin{align*}
    &ds^2_t
        =(\cosh^2 t) ds^2+ 2(\cosh t\sinh t) dh^2+
          (\sinh^2 t) \langle d\nu,d\nu \rangle, \\
    &dh^2_t
        =(\cosh t \sinh t)ds^2+2(\cosh^2 t+\sinh^2 t)dh^2+
          (\cosh t \sinh t)\langle d\nu,d\nu \rangle, \\
    &\langle d\nu_t,d\nu_t \rangle
        =
          (\sinh^2 t) ds^2+ 2(\cosh t\sinh t) dh^2+
          (\cosh^2 t) \langle d\nu,d\nu \rangle, 
\end{align*}
where $\langle d\nu,d\nu \rangle$ is the third
fundamental form of $f$.
Since $d\sigma^2=ds^2+\langle d\nu,d\nu \rangle$, we have
\[
   dh^2_t
     =(\cosh t \sinh t)d\sigma^2+2(\cosh^2 t+\sinh^2 t)dh^2, 
\]
and
\[
    d\sigma^2_t:=
        ds^2_t + \langle d\nu_t,d\nu_t \rangle
     =(\cosh^2 t+\sinh^2 t)d\sigma^2+4\cosh t\sinh t \, dh^2.
\]
Since $dh^2$ is proportional to $d\sigma^2$,
$dh^2_t$ and $d\sigma^2_t$ are also proportional.
Moreover, since $f$ is not a horosphere,
$ds^2$ is not proportional to $dh^2$ and thus
$f_t$ is flat for all $t\in \R$.
\end{proof}

From now on, we assume that $M^2$ is oriented.
(If $M^2$ is not orientable, we can take the double cover.)
Then there is a complex structure on $M^2$ such that $d\sigma^2$ is
hermitian.
Since the second fundamental form is proportional to $d\sigma^2$,
this complex structure of $M^2$ coincides with the one treated in
Section~\ref{sec:preliminaries}, as long as $f$ is an immersion.
So, we shall call this complex structure 
the {\it canonical complex structure\/},
and $M^2$ is always considered as a Riemann surface.
Now the following assertion holds.

\begin{proposition}
\label{prop:2-equiv}
 Let $M^2$ be a Riemann surface and 
 $E\colon \widetilde M^2\to \SL(2,\C)$ be a holomorphic Legendrian
 immersion defined on the universal cover $\widetilde M^2$ such that
 $f=EE^*$ is single-valued on $M^2$.
 Then $f$ is a flat front. %
 Moreover, if we set
 \begin{equation}\label{eq:diffE}
    E^{-1}dE=
        \begin{pmatrix}
	 0 & \theta \\
	 \omega & 0
	\end{pmatrix},
 \end{equation}
 the first and the second fundamenal forms are represented as
\begin{equation}\label{eq:front-fundforms}
 \begin{aligned}
   ds^2 &= \omega\theta+\bar\omega\bar\theta+(|\omega|^2+|\theta|^2),\\
   dh^2 &= |\theta|^2-|\omega|^2.
 \end{aligned}
\end{equation}
 Conversely, any flat front is given as a projection of a holomorphic
 Legendrian immersion.
\end{proposition}

\begin{proof}
Let 
\[
    E\colon M^2\longrightarrow \SL(2,\C)
\]
be a holomorphic Legendrian map and 
$(\omega,\theta)$ as in \eqref{eq:diffE}.
Then $E$ is an immersion if $|\omega|^2+|\theta|^2$ is positive
definite.
On the other hand, we have
\begin{align*}
 df
   &=dE\,E^*+E\,dE^*
     =E(E^{-1}dE + (E^{-1}dE)^*)E^*\\
     &= E\begin{pmatrix}
	  0 & \theta+\bar\omega \\
          \omega+\bar\theta & 0 
        \end{pmatrix}E^*,\\
 d\nu&=
     dEv_3E^*+Ev_3dE^*=
         E(E^{-1}dEv_3 + v_3(E^{-1}dE)^*)E^*\\
     &=  E\begin{pmatrix}
	  0 & -\theta+\bar\omega \\
          \omega-\bar\theta & 0 
        \end{pmatrix}E^*.
\end{align*}
In the identification as in  \eqref{eq:lor-herm}, 
the canonical Lorentzian inner product
is given as follows
\[
  \langle X,Y\rangle:=-\frac12\trace(X\widetilde Y)\qquad 
       \bigl(X,Y\in \Herm(2)\bigr),
\]
where $\widetilde Y$ is the cofactor matrix of $Y$.
If $Y\in \SL(2,\C)$, we have $\tilde Y=Y^{-1}$.
 Then we have
\begin{align*}
   ds^2 &= \langle df,df\rangle
        =-\frac{1}{2}\trace\left\{
          \begin{pmatrix}
	   0 & \theta+\bar\omega \\
           \omega+\bar\theta & 0 
	  \end{pmatrix}
          \begin{pmatrix}
	   0 & -\theta-\bar\omega \\
           -\omega-\bar\theta & 0 
	  \end{pmatrix}\right\}\\
       &=(\omega+\bar\theta)(\bar\omega+\theta).
\end{align*}
 Similarly, since $dh^2=-\langle df,d\nu\rangle$, we have
 \eqref{eq:front-fundforms}.
Thus,
the pull-back of the Sasakian metric by $(f,\nu)$ as in 
\eqref{eq:sasakian} is represented as
\begin{equation}\label{eq:sasaki-holo}
 d\sigma^2 = \langle df,df\rangle + \langle d\nu,d\nu\rangle
           = 2(|\omega|^2+|\theta|^2).
\end{equation}
Hence $L_f$ is an immersion if and only if $|\omega|^2+|\theta|^2$
is positive definite. This proves the assertion.
\end{proof}
\begin{remark}
 As seen in the proof of Proposition~\ref{prop:2-equiv},
 the $(1,1)$-part of the first fundamental form
\begin{equation}\label{eq:one-one}
 ds^2_{1,1}=|\omega|^2+|\theta|^2
\end{equation}
 is equal to the one half of $d\sigma^2$, the pull-back of 
 the Sasakian metric by $(f,\nu)$.
 Moreover, $ds^2_{1,1}$ is the pull-back of the bi-invariant
 Hermitian metric of $\SL(2,\C)$ by $E$.
\end{remark}
 
Now, a flat front $f\colon M^2\to H^3$ can be interpreted from
two different points of view.
The first is the projection of a (real) Legendrian immersion
$L_f\colon M^2 \to T_1^* H^3\cong \F$,
and the second is the projection of a (holomorphic)
Legendrian immersion $E_f\colon \widetilde M^2\to \SL(2,\C)$. 
One can naturally expect that these two Legendrian
immersions are related. In fact, $\SL(2,\C)$ acts $\F$
transitively and we can write 
\[
   T_1^*H^3\cong T_1H^3\cong \mathcal F
         \cong \SL(2,\C)/\U(1).
\]
We denote the canonical projection by
\[
   p^{}_{\SL}\colon \SL(2,\C)\longrightarrow \F.
\]
We prove the following: 

\begin{proposition}\label{prop:contact}
 The pull-back of the contact form $\eta$ by  $p^{}_{\SL}$ is equal to
 the real part of the holomorphic contact form on $\SL(2,\C)$,
 that is, 
 \[
    p^{*}_{\SL}(\eta)
      =2\Re(s_{22}ds_{11}-s_{12}ds_{21}) 
 \]
 holds, where $(s_{ij})\in \SL(2,\C)$.
 In particular, the real Legendrian immersion $L_f$ can be
 interpreted as the projection of  a holomorphic Legendrian immersion
 $E_f$.
\end{proposition}

\begin{proof}
Since
\[
   x^{-1}=(ss^*)^{-1}=(s^*)^{-1}s^{-1},
     \qquad dv=s(s^{-1}dsv_3+v_3(s^{-1}ds)^*)s^*,
\]
we have
\begin{align*}
  \eta&=\langle x,dv \rangle=-\frac12\trace(x^{-1}dv)\\
      &=-\frac12\trace(s^{-1}ds\, v_3+v_3(s^{-1}ds)^*)\\
      &=-\Re\left[\trace(s^{-1}ds\,v_3)\right]\\
      &=2\Re(s_{22}ds_{11}-s_{21}ds_{12}),
\end{align*}
where we set $s=(s_{ij})\in \SL(2,\C)$. 
\end{proof}
\begin{remark}
 Since the holomorphic Legendrian lift $E_f$ of a flat front $f$ is not
 single-valued on $M^2$ in general, $E_f$ has  
 the monodromy representation
 $\rho_{f} \colon \pi_1(M^2)\to\SU(2)$ such that 
 $E_f\circ\tau=E_f\rho_f(\tau)$ for any deck transformation
 $\tau\in\pi_1(M)$.
 On the other hand, since $E_f\circ\tau$ is also Legendrian,
 the representation $\rho_f$ is reducible, that is, it
 reduces to the isotropy group $\U(1)$ of the action
 of $\SL(2,\C)$ to $\F$.
\end{remark}

We can define the hyperbolic Gauss map of the flat front $f$ just 
the same way as for an immersion, that is, 
\[
   G=\frac{A}{C},\qquad
   G_*=\frac{B}{D},
\]
where 
\[
  E_f=
  \begin{pmatrix} 
   A & B \\ C & D 
  \end{pmatrix}.
\]
These are single-valued on $M^2$.
Since $AD-BC=1$, $G(p)$ never coincides with $G_*(p)$,
and we get the holomorphic map
\[
  \G:=(G,G_*)\colon
     M^2 \longrightarrow S^2\times S^2\setminus
      \{\text{the diagonal set}\}=:\Geod(H^3),
\]
where $\Geod(H^3)$ is the set of oriented geodesics in $H^3$.

Now, we shall prove the following:
\begin{theorem}\label{thm:imm}
 Let $\widetilde M^2$ be the universal cover of a Riemann surface 
 $M^2$ and $E\colon \widetilde M^2 \to \SL(2,\C)$ a holomorphic
 Legendrian map such that $f=EE^*$ is single-valued on $M^2$. 
 Then the following four assertions are mutually equivalent{\rm :}
 \begin{enumerate}
  \item\label{item:imm-1} 
       $E$ is an immersion.
  \item\label{item:imm-2}
       $L_f$ is an immersion.
  \item\label{item:imm-3}
       The $(1,1)$-part of the first fundamental form
       \[
        ds^2_{1,1}=
        \left| \omega \right|^2+\left| \theta \right |^2
       \]
       is positive definite, where $\omega$ and $\theta$ are
       the off-diagonal components of $E^{-1}dE$.
 \item\label{item:imm-4}  
       $\G:=(G,G_*)\colon M^2 \to \Geod(H^3)$ is an immersion.
 \end{enumerate}
\end{theorem}
\begin{remark}\label{rem:branch}
 A point  where  $E$ degenerates is called a {\it branch point\/}
 of $E$ (or $f=EE^*$), and the projection of a
 holomorphic Legendrian curve is called a {\it branched flat front\/}. 
 The above conditions imply that $f$ is free of branch point.
\end{remark}

\begin{proof}[Proof of Theorem~\ref{thm:imm}]
The equivalency of the first three assertions follows from the proof of
 Proposition~\ref{prop:2-equiv}.
 So it is sufficient to prove that \ref{item:imm-3} and \ref{item:imm-4}
 are equivalent.
 By \eqref{eq:omega} and \eqref{eq:theta}, we have the following identity
 \[
    ds^2_{1,1}=|\omega|^2+|\theta|^2
      =|C^2dG|^2+|D^2dG_*|^2.
 \]
 If both $C$ and $D$ are non-vanishing, then the equivalency of 
 \ref{item:imm-3} and \ref{item:imm-4} are obvious. 
 So we consider the case $C=0$ (resp.\ $D=0$).
 In this case $D\ne 0$ (resp.\ $C\ne 0$).
 Then $\omega=0$ (resp.\ $\theta=0$) and
 \[
    ds^2_{1,1}=|\theta|^2=|D^2dG_*|^2,
      \qquad (\mbox{resp. }
      ds^2_{1,1}=|\omega|^2=|C^2dG|^2).
 \]
 Then $ds^2$ is positive definite
 if and only if $dG_* \ne 0$ (resp.\ $dG \ne 0$).
 This proves the assertion.
\end{proof}

In \cite{KUY}, the authors gave a representation formula for Legendrian
curves in $\SL(2,\C)$ via the data $(G,G_*)$. 
We reformulate it for construction of flat fronts in $H^3$, as follows.

\begin{theorem}\label{thm:rep}
 Let $G$ and $G_*$ be non-constant meromorphic functions on a Riemann
 surface $M^2$ such that $G(p)\ne G_*(p)$ for all $p\in M^2$.
 Assume that
 \begin{enumerate}
  \item\label{item:rep-1}
      all poles of the $1$-form $\dfrac{dG}{G-G_*}$ are of order $1$,
      and
  \item\label{item:rep-2}
      $\displaystyle\int_{\gamma}\dfrac{dG}{G-G_*}\in i \R$ 
      holds for each loop $\gamma$ on $M^2$.
 \end{enumerate}
 Set
 \begin{equation}\label{eq:xi}
    \xi(z) := c \exp\int_{z_0}^{z}\frac{dG}{G-G_*},
 \end{equation}
 where $z_0\in M^2$ is a base point and $c\in \C\setminus\{0\}$ is
 an arbitrary constant.
 Then 
 \begin{equation}\label{eq:repG2}
     E := \begin{pmatrix}
      G/\xi &      \xi G_{*}/(G-G_{*}) \\
      1/\xi &      \hphantom{G_*}\xi/(G-G_{*}) 
     \end{pmatrix}
 \end{equation}
 is a non-constant meromorphic Legendrian curve  
 defined on $\widetilde M^2$ in $\PSL(2,\C)$ whose hyperbolic Gauss maps
 are $G$ and $G_*$, such that the projection $f=EE^*$  is single-valued on
 $M^2$. 
 Moreover, $f$ is a front if and only if $G$ and $G_*$ have no common
 branch points.  
 Conversely, any non-totally-umbilical flat fronts can be constructed in
 this manner.
\end{theorem}

\begin{proof}
If we give a pair $(G,G_*)$ of non-constant meromorphic functions,
on a Riemann surface $M^2$ satisfying \ref{item:rep-1},
a meromorphic map $E$ defined by \eqref{eq:repG2} is a holomorphic
Legendrian curve in $\SL(2,\C)$,
which is a consequence of Theorem 3 of \cite{KUY}.
Then the second condition (2) implies that $f=EE^*$ is single-valued on
$M^2$.
Now, by Theorem~\ref{thm:imm}, the branched flat front $f$ is free of
branch point if and only if the pair $(G,G_*)$ gives an immersion of
$M^2$ into $S^2\times S^2$.

Since any flat front can be lifted to a holomorphic Legendrian curve
defined on $\widetilde M^2$, Theorem~3 of \cite{KUY} yields also
that any non-totally-umbilical flat front can be constructed
in this manner. 
(If one of the hyperbolic Gauss maps is constant, 
 it is totally umbilic, i.e.\ locally a horosphere.)
\end{proof}

\section{Flat fronts with complete ends}
\label{sec:complete}
We define completeness of fronts as follows: 
\begin{definition}
 Let $M^2$ be a $2$-manifold.
 A front $f\colon M^2 \to H^3$ is called {\it complete\/} if
 there is a symmetric 2-tensor $T$ with compact support
 such that the sum 
 \[
   T+ds^2
 \]
 is a complete Riemannian metric of $M^2$, 
 where $ds^2$ is the first fundamental form of $f$.
\end{definition}
\begin{remark}
 It should be noted that the parallel family of a given complete
 front $f$ may contain an incomplete flat front. 
 For example, the hyperbolic cylinder, 
 that is the surface equidistance from a geodesic 
 (see Example~\ref{ex:ffofrev} in Section~\ref{sec:example}),
 contains a geodesic in its parallel family.
\end{remark}
The following assertion is a simple consequence of Lemma 2 of
\cite{GMM}.

\begin{lemma}\label{lem:finite-topology}
 Let $M^2$ be an oriented $2$-manifold and
 $f\colon M^2 \to H^3$ a complete flat front.
 Then there exists a compact Riemann surface
 $\overline M^2$ and finitely many points 
 $p_1,\dots,p_n \in  \overline{M}^2$ such that $M^2$
 {\rm (}as a Riemann surface{\rm )} is biholomorphic
 to $\overline M^2\setminus \{p_1,\dots,p_n\}$.
 Moreover, the Hopf differential $Q$ of $f$
 can be extended meromorphically on $\overline M^2$.
\end{lemma}

These points $p_1,\dots,p_n$ are called {\it ends\/} of the front $f$.
\begin{proof}[Proof of Lemma~\ref{lem:finite-topology}]
Since $f$ is complete, there exists a symmetric tensor $T$ whose support
is a compact subset of $M^2$ such that
\[
    d\bar s^2:=T+ds^2
\]
is complete, where $ds^2$ is the first fundamental form. 
Since the Gaussian curvature of $ds^2$ vanishes, the total absolute
curvature of $d\bar s^2$ is finite. 
Then by Huber's theorem, there is a compact 2-manifold $\overline M^2$
and finite points $p_1,\dots,p_n\in\overline{M}^2$ such that $M^2$ is
diffeomorphic to $\overline M^2\setminus \{p_1,\dots,p_n\}$.
Now we take a sufficiently small neighborhood $U_j$ of an end $p_j$
such that  $d\bar s^2=ds^2$ holds on $U_j$.
If $|\omega|=|\theta|$ holds on a point $q\in U_j$, $ds^2$ degenerates
at $q$ because of \eqref{eq:first}.
Hence $|\omega|\neq|\theta|$ holds on $U_j$.
If $|\omega|>|\theta|$, by (25) of \cite{GMM}, we have
\[
   ds^2 = \omega\theta+\bar\omega\bar\theta+|\omega|^2
                +|\theta|^2 
        \leq 2|\omega||\theta|+|\omega|^2
                +|\theta|^2 
        =(|\omega|+|\theta|)^2 \leq  4|\omega|^2.
\]
Since $ds^2$ is complete at $p_j$, so is a metric $|\omega|^2$.
Moreover, by a holomorphicity of $\omega$ with respect to the 
complex structure induced from the second fundamental form, 
the metric $|\omega|^2$ is a flat metric which is conformal
to the complex structure of $M^2$ and complete at $p_j$.
This proves the first assertion.
In the case of $|\omega|<|\theta|$, we have the conclusion
using a conformal metric $|\theta|^2$.

The meromorphicity of $Q$ is proved in Lemma~2 of \cite{GMM}.
\end{proof}
As seen in the proof of Lemma~\ref{lem:finite-topology}, we have
\begin{corollary}\label{cor:complete-sasaki}
 If a flat front $f$ is complete, so is the $(1,1)$-part 
\[
    ds^2_{1,1}=|\omega|^2+|\theta|^2
\]
 of the first fundamental form.
\end{corollary}

\begin{figure}
\begin{center}
 \includegraphics[width=3cm]{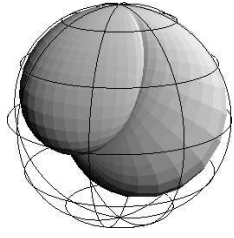}
\end{center}
\caption{An incomplete flat front (Remark~\ref{rem:incomplete}).}
\label{fig:incomplete}
\end{figure}
\begin{remark}\label{rem:incomplete}
 If a flat front is complete and is also a proper mapping, 
 then its image is a closed subset of $H^3$. 
 However, a proper flat front $f$ whose image is  closed in $H^3$
 may not be complete.
 (When  $f$ has no singularity, it is complete by the Hopf-Rinow
 theorem.)
 In fact, we consider a flat front $f=EE^*:\C\to H^3$, where
 \begin{equation}\label{eq:incomplete}
   E:=\begin{pmatrix}
       z e^{-z} & (z-1) e^z \\
       \hphantom{z}e^{-z} & \hphantom{(z-1)}e^z
      \end{pmatrix}.
 \end{equation}
 It can be easily checked that $f(z)$ tends to the north pole of 
 the ideal boundary in the Poincare ball as $z\to \infty$,
 which implies that $f$ is a proper mapping.
 But the first fundamental form vanishes on the
 imaginary axis, which appears as cuspidal edges.
 (See Figure~\ref{fig:incomplete}. 
  The criterion for singularities of flat fronts will be
  appeared in the forthcoming paper \cite{KRSUY}.)
\end{remark}
Although it is obvious that 
there are no compact flat surfaces in $H^3$
according to the classical fact, we can also prove the 
non-existence of compact flat fronts. 
\begin{proposition}\label{prop:noncompact}
  There are no compact flat fronts without boundary. 
\end{proposition}
\begin{proof}
  Suppose that $f\colon M^2 \to H^3$ is a compact flat front.
  Take a holomorphic Legendrian lift $E$ of $f$.
  Let us be conscious that $f$ and $E$ are matrix-valued. 
  Then the trace of $f$ satisfies 
 \[
   (\trace f)_{z\zb}=\trace (f_{z\zb})=
      \trace\left\{E_z(E_z)^*\right\}  
      \ge 0, 
 \]
  where $z$ is a complex coordinate of $M^2$.
  Hence the function $\trace f\colon{}M^2\to\R$ is subharmonic,
  which must be constant, since $M^2$ is compact. 
  By an isometry in $H^3$, we may assume that $f(z_0)=I$, where 
  $I$ is the $2 \times 2$ identity matrix. 
  Then $\trace f$ is identically $2$. 
  On the other hand, $\det f$ is identically $1$. 
  This implies that 
  the eigenvalues $\lambda_1, \lambda_2$ of $f$ satisfy  
 \[
    \lambda_1+\lambda_2 =2, \qquad \lambda_1 \lambda_2 =1. 
 \]
 Hence we have $\lambda_1=\lambda_2=1$.
  Since $f$ is Hermitian, this implies that $f(z)$ is equal to the
  identity matrix, a contradiction.
\end{proof}

G\'alvez, Mart\'\i{}nez and Mil\'an investigated complete ends of flat
surfaces deeply.
The following fact is proved in \cite{GMM}.

\begin{lemma}[Theorem 4 of \cite{GMM}]
\label{lemma:GMM}
 Let $p$ be an end of a complete flat front.
 Then the following three conditions are mutually
 equivalent.
\begin{enumerate}
 \item\label{item:GMM-1}
      The Hopf differential $Q$ has at most a pole of order $2$ at $p$.
 \item\label{item:GMM-2}
      The positive hyperbolic Gauss map $G$ has
      at most a pole at $p$.
 \item\label{item:GMM-3} 
      The negative hyperbolic Gauss map $G_*$ has at most a pole at $p$.
\end{enumerate}
\end{lemma}
\begin{remark}
 The hyperbolic Gauss maps and Hopf differential of the flat front 
 $f=EE^*$ as in \eqref{eq:incomplete} are
 \[
    G=z,\qquad G_*=z-1,\qquad Q=dz^2.
 \]
 This means that meromorphicity of $G$ and $G_*$ does not imply that $Q$
 has at most poles of order $2$ without assuming the completeness of
 ends.
 In fact, $Q$ has pole of order $4$ at $z=\infty$.
\end{remark}

If an end of a flat front satisfies one of the above three conditions,
it is called a {\it regular\/} end. 
If an end is not regular, it is called an {\it irregular\/} end.
An end $p$ is said to be {\it embedded\/} if there is a 
neighborhood $U$ of $p\in\overline M^2$, 
such that the restriction of the front to $U\setminus \{p\}$
is an embedding.

The following fact is also in \cite{GMM}.

\begin{lemma}[Theorem 5 of \cite{GMM}]
\label{embedded_1} 
 Let $p$ be a regular end of a complete flat front.
 Suppose that $|\theta|<|\omega|$ at $p$.
 Then $p$ is embedded if and only if
 it is not a branch point of the positive hyperbolic 
 Gauss map $G$.
\end{lemma}

We also need the following:
\begin{lemma}\label{lemma:same}
 Let $p$ be a regular end of complete flat front,
 then $G(p)=G_*(p)$ holds.
 Namely, the two Gauss maps take the same value at a regular
 end. 
\end{lemma}

\begin{proof}
 Assume that $G(p)\neq G_*(p)$ for a regular end $p$.
 By Lemma~\ref{lemma:GMM}, 
 $G(z)$ and $G_*(z)$ are both meromorphic at $p$.
 In particular, the function $\xi(z)$ defined in 
 Theorem \ref{thm:rep} is holomorphic. 
 Then so is $E$ and this contradicts the completeness of the first
 fundamental form of the front at $p$.
 Thus we have $G(p)=G_*(p)$.
\end{proof}

Let $f\colon \overline M^2\setminus \{p_1,\dots,p_n\}\to H^3$ be 
a complete flat front whose first and second canonical forms are
$\omega$ and $\theta$ respectively.
Suppose that all ends $p_1,\dots,p_n$ are regular. 
By Lemma~2 of \cite{GMM}, there exist real numbers $\mu_j^{}$
and $\mu_j^*$ ($j=1,\dots,n$)
such that 
\begin{alignat*}{2}
 \omega(z)&=(z-p_j)^{\mu_j}\omega_0(z) \qquad   &(\omega_0(p_j)&\ne 0),\\
 \theta(z)&=(z-p_j)^{\mu_j^*}\theta_0(z) \qquad &(\theta_0(p_j)&\ne 0),
\end{alignat*}
where 
$(U;z)$ is a complex coordinate around $p_j$, and
$\omega_0(z)$ and $\theta_0(z)$ are holomorphic 1-forms
defined on $U$.
The real number $\mu_j$ (resp.\ $\mu_j^*$) does not depend on the
choice of the coordinate $z$ and is equal to the order of the
pseudometric $|\omega|^2$ (resp.\ $|\theta|^2$),
namely
\begin{equation}\label{eq:omega-theta-expansion}
\begin{aligned}
   |\omega|^2&=|z-p_j|^{2\mu_j^{}}\left(a_j^{}+o(1)\right)|dz|^2,\\
   |\theta|^2&=|z-p_j|^{2\mu_j^*}\left(a_j^*+o(1)\right)|dz|^2,
\end{aligned}
\end{equation}
where $a_j^{}$ and $a_j^*$ are positive real numbers and 
$o(1)$ denotes higher order terms.
By \eqref{eq:Hopf}, we have
\begin{equation}\label{eq:Hopf-order}
 \mu_j^{}+\mu_j^* = \ord_{p_j} Q,
\end{equation}
where $\ord_{p_j}Q$ is the order of the Hopf differential
$Q$ at $p_j$.
Suppose that the Hopf differential $Q$ has the following Laurent
expansion at $z=p_j$:
\[
   Q=\frac{1}{(z-p_j)^2}\bigl(q_{-2}(p_j)+o(1)\bigr)dz^2,
\]
where $o(1)$ is a function satisfying 
$\displaystyle \lim_{z\to p_j}o(1)=0$.
The following lemma is a direct consequence of
the formula \eqref{eq_sch}.
\begin{lemma}\label{lem_relation}
 The identity
 \[
        4q_{-2}(p_j)
          =m_j^{}(m_j^{}+2)-\mu_j^{}(\mu_j^{}+2)
          =m_j^{*}(m_j^{*}+2)-\mu_j^{*}(\mu_j^{*}+2)
 \]
 holds,
 where $\mu_j^{}$  {\rm (}resp.\ $\mu_j^*${\rm)} is
 the oreder of the pseudometric $|\omega|^2$ 
 {\rm(}resp.\ $|\theta|^2${\rm)}, 
 and $m_j$ {\rm (}resp.\ $m_j^*${\rm )} is the branching order of $G$
 {\rm (}resp.\ $G_*${\rm )} at $p_j$, respectively. 
 {\rm (}e.g.\ $m_j=1$ implies $p_j$ is a double point of $G$.{\rm )}
\end{lemma}

Now we shall prove the following:
\begin{proposition} \label{embedded_2}
 A regular end $p_j$ of a complete flat front is embedded if and only if 
 either $G$ or $G_*$ does not branch at $p_j$.
\end{proposition}
\begin{proof}
 First, we assume the end $p_j$ is embedded.
 As pointed out in previous paper \cite{KUY}, the holomorphic Legendrian
 lift $E_f$ has the following duality:  We set
 \[
     \hat E_f:=
          E_f  \begin{pmatrix}
		0 & i \\
		i & 0
	       \end{pmatrix}.
 \]
 Then $\hat E_f$ is also a holomorphic Legendrian immersion
 such that $f=\hat E_f^{}\hat E_f^*$ and the role of $(G,\omega)$ 
 and $(G_*, \theta)$ interchanges.
 Then, by replacing $E_f$ by $\hat E_f$ if necessary,
 we may assume $|\theta|<|\omega|$ near $p_j$.
 Then by Lemma~\ref{embedded_1}, $G$ does not branch at $p_j$.

 Conversely, we assume either $G$ or $G_*$ does not branch at $p_j$. 
 By replacing $E_f$ by $\hat E_f$ if necessary,
 we may assume $G$ does not branch at $p_j$, that is, $m_j=0$ holds.
 If $|\theta|<|\omega|$ near $p_j$, the assertion follows from
 Lemma~\ref{embedded_1}.
 So we may assume $|\theta|>|\omega|$ near $p_j$,  and then we
 have $\mu_j^{}\ge \mu_j^*$.
 By Lemma~\ref{lem_relation}, we have
 \[
     (m_j^{}-m_j^*)
     (m_j^{}+m_j^*+2) =
     (\mu_j^{}-\mu_j^*)
     (\mu_j^{}+\mu_j^*+2) .
 \]
 By \eqref{eq:Hopf-order}, we have
 \[
       \mu_j+\mu^*_j+2=\ord_{p_j}Q+2.
 \]
 Since $p_j$ is regular, $\ord_{p_j}Q\ge -2$. 
 Thus we have $0=m_j\ge m^*_j(\ge 0)$ and so both of $G$ and $G_*$
 do not branch at $p$.
 So by replacing $E_f^{}$ by $\hat E_f$ if necessary,
 we may assume $|\theta|<|\omega|$ and get the embeddedness
 of $p$ directly from Lemma~\ref{lem_relation}.
\end{proof}

Now we shall prove the following assertion, stated in the introduction.

\begin{theorem}\label{thm:Oss-ineq}
 Let
 $f\colon \overline M^2\setminus \{p_1,\dots,p_n\}\to H^3$ 
 be a complete flat front whose ends are all regular.
 Then the following inequality holds
 \[
   \deg G+\deg G_*\ge n,
 \]
 and equality holds if and only if all ends are embedded.
\end{theorem}

To prove the theorem, we shall prepare two lemmas.
\begin{lemma} \label{lemma-new1}
 Let $g$ and $h$ be meromorphic functions on a compact Riemann
 surface $\overline M^2$.
 Suppose that $g$ and $h$ have no common poles, then
 \[
    \deg(ag+bh)=\deg g+\deg h
 \]
 unless $ab=0$, where $a$, $b\in\C$ are constants.
\end{lemma}
\begin{proof}
 Since $g,h$ are meromorphic, each of their degrees
 is equal to the number of their poles counting multiplicities.
 If we denote the divisor of poles of $g$ and $h$
 by 
 \[
   P(g):=s_1q_1+\cdots +s_n q_n, \qquad
   P(h):=l_1r_1+\cdots +l_k r_k,
 \]
 where $\{q_1,\dots,q_n\}$ (resp.\ $\{r_1,\dots,r_k\}$)
 is the set of poles of $g$ (resp.\ $h$) and $s_j$, $l_i$ are 
 positive integers.
 Then we have $P(a g+bh)=P(g)+P(h)$ unless $ab=0$.
 Thus
 \[
   \deg g+\deg h=\sum_{j=1}^n s_j + \sum_{i=1}^k l_i
     =\deg(ag+bh).
 \]
 This proves the lemma.
\end{proof}
\begin{lemma} \label{lemma-new2}
 Let  $f\colon \overline M^2\setminus \{p_1,\dots,p_n\}\to H^3$
 be a complete flat front.
 Suppose that an end $p=p_j$ is regular.
 Then $p$ is an embedded end if and only if
 the difference of two Gauss maps
 \[
     h:=G-G_*
 \]
 does not branch at $p$.
\end{lemma}
\begin{proof}
 If $p$ is not a branch point of $h$, 
 then either $G$ or $G_*$ does not branch at $p$.
 Then embeddedness of the end $p=p_j$ follows from
 Proposition~\ref{embedded_2}.
 Conversely, suppose now that an end $p$ is embedded.
 We take a complex coordinate $z$ around $p$ such that $z(p)=0$.
 Then, by an isometry of $H^3$, we may assume
 $G(0)=G_*(0)=0$ because  of Lemma~\ref{lemma:same} and
 \eqref{eq:gauss-change}.
 It follows from Proposition~\ref{embedded_2} that $G$ and $G_*$
 are expanded as
 \[
   G(z)= a\,z+o(z)\qquad
    \text{and}\qquad
   G_*(z) = a_*z+o(z),
 \]
 where 
 $a$ and $a_*$ are complex numbers such that $a\neq 0$
 or $a_* \neq 0 $ holds, and $o(z)$ denotes a higher order term.
 Thus by \eqref{eq:Hopf}, the Hopf differential $Q$ is expanded as
 \[
   Q = -\frac{aa_*+o(1)}{\{(a-a_*)z+o(z)\}^2}\,dz^2.
 \]
 Hence by Lemma~\ref{lemma:GMM}, we have $aa_*=0$ or 
 $a-a_*\neq 0$.
 If $a-a_*\neq 0$, it follows that $h$ does not branch 
 at $0$.
 In the case of  $aa_*=0$, one of $G$ and $G_*$ branches at $0$ and the
 other does not.
 Then $h=G-G_*$ does not branch at $0$.
\end{proof}
\begin{proof}[Proof of Theorem~\ref{thm:Oss-ineq}]
Taking an isometry if necessary, 
 we may assume that all ends $p_1,\dots,p_n$ are not poles of both of two
 Gauss maps $G$ and $G_*$.
 Since the ends of the front are equal to the 
 zeros of $h:=G-G_*$, $G$ and $G_*$ has no common poles.
 Moreover, the zero divisor $Z(h)$ of
 the meromorphic function $h$ is of the form
 \[
     Z(h)=\sum_{j=1}^n m_jp_j, 
 \]
 where $m_1,\dots,m_n$ are positive integers.
 Then by Lemma \ref{lemma-new1}
 \[
   \deg G+\deg G_*=\deg h
         =\sum_{j=1}^n m_j\ge n,
 \]
 which proves the inequality.
 Moreover, the equality holds if and only if 
 \[
    m_1=\dots=m_n=1.
 \]
 This holds if and only if all ends are embedded by 
 Lemma \ref{lemma-new2}.
\end{proof}

\begin{remark}\label{rem:osserman-branched}
 As seen in the proof,  the inequality of 
 Theorem~\ref{thm:Oss-ineq} holds even if 
 $f$ has branch points (see Remark~\ref{rem:branch}), 
 that is, common branch points  of $G$ and $G_*$ on $M$.
 However, the category of branched flat front seems too wide for the
 study of flat surfaces.
 In fact, branched covers of flat fronts are all
 branched flat fronts whose images are
 the same as the original fronts.
\end{remark}
\begin{remark}
 Let $\overline{M}^2$ be a compact Riemann surface  with positive genus.
 Since there are no meromorphic functions on $\overline{M}^2$
 of degree $1$,
 a complete flat front defined on $\overline{M}^2$ excluded
 finite number of points must have at least $4$ ends.
 It is interesting to consider a problem 
 `{\it Is there a flat front with positive genus with
          $4$ embedded ends?}' %
 Though there exists a genus one flat front with $5$ ends
 as shown in Example~\ref{ex:torus} later,
 it is unknown whether there exists a genus $1$ flat front 
 with $4$ ends.
 On the other hand, one can construct a genus $1$ 
 {\em branched\/} flat front with $4$ embedded ends, 
 but the image of such a front is a double cover of
 a genus $0$ flat  front.
\end{remark}

As seen in Remark~\ref{rem:osserman-branched},
the inequality of Theorem~\ref{thm:Oss-ineq} is valid for 
branched flat fronts.
On the other hand, we show an inequality which reflects properties 
of fronts.
Before stating the result, we prepare a terminology:
\begin{definition}
 Let $p$ be a regular end of a flat front and $(\omega,\theta)$ as in
 \eqref{eq:diffE}.
 Then the end $p$ is called {\em cylindrical\/} if 
\[
    \ord_p|\omega|^2=\ord_p|\theta|^2=-1
\]
 holds, where $\ord_p|\omega|^2$ \rm{(}resp.\ $\ord_p|\theta|^2${\rm)}
 are the order of a pseudometric $|\omega|^2$ {\rm(}resp.\ $|\theta|^2$)
 as in \eqref{eq:omega-theta-expansion}.
\end{definition}
The ends of a hyperbolic cylinder, a surface equidistance from a
geodesic, are cylindrical.
(See Example~\ref{ex:ffofrev}. 
 see also \cite[page~427]{GMM}, or \cite[Example~4.1]{KUY}.)
By Theorem~6 of \cite{GMM}, we have
\begin{lemma}[{\cite[Theorem 6]{GMM}}]%
\label{lem:cylindrical-assymptotic}
 A cylindrical end is asymptotic a finite cover of a hyperbolic 
 cylinder.
\end{lemma}
\begin{lemma}\label{lem:cylindrical}
 A regular end $p$ of a complete flat front is cylindrical if 
 and only if $\ord_p|\omega|^2=-1$ or $\ord_p|\theta|^2=-1$ holds.
\end{lemma}
\begin{proof}
 Assume $\ord_p|\omega|=-1$.
 Then $\omega$ is written as 
 $\omega = (z-p)^{-1}\omega_0(z)$, where $\omega_0$ is a holomorphic
 $1$-form such that $\omega_0(p)\neq 0$ and $z$ is a complex coordinate 
 around $p$.
 On the other hand, by isometry of $H^3$ if necessary, we 
 may assume $G(p)\neq 0$ because \eqref{eq:gauss-change}.
 Then $G$ is is written as $(z-p)^mG_0(z)$, where $m\geq 1$ is an
 integer and $G_0(z)$ is a holomorphic function such that $G_0(p)\neq 0$.
 Then by \eqref{eq_sch}, we conclude that the order of the Hopf
 differential $Q$ at $z=p$ is $-2$.
 Hence by \eqref{def_hopf}, $\ord_p|\theta|^2=-1$ holds.
 Similarly, if $\ord_p|\theta|^2=-1$, we have $\ord_p|\omega|^2=-1$.
\end{proof}
Let $f\colon{}\overline M^2\setminus\{p_1,\dots,p_n\}\to H^3$ be a 
complete flat front and $(\omega,\theta)$ as in \eqref{eq:diffE}.
Then, by \eqref{eq:sasaki-holo},
the pull-back of the Sasakian metric of $T_1^*H^3$ by the Legendrian
lift of $f$ is
\begin{equation}\label{eq:sasaki-metric}
   d\sigma^2= 2(|\omega|^2+|\theta|^2).
\end{equation}
Since $\omega$ and $\theta$ are holomorphic $1$-forms,  the 
Gaussian curvature $K_{d\sigma^2}$ of $d\sigma^2$ is non-positive.
Moreover, $d\sigma^2$ is complete because of
Corollary~\ref{cor:complete-sasaki}.
\begin{proposition}\label{prop:osserman-curvature}
 Let $f\colon{}M=\overline M^2\setminus\{p_1,\dots,p_n\}\to H^3$ be
 a complete flat front and $d\sigma^2$ the pull-back of the Sasakian 
 metric of $T_1^*H^3$ as in \eqref{eq:sasaki-metric}.
 Then total Gaussian curvature of $d\sigma^2$ is finite and 
\begin{equation}\label{eq:osserman-curvature}
   \frac{1}{2\pi}\int_{M^2}(-K_{d\sigma^2})dA_{d\sigma^2}
       \geq -\chi(M^2) = -\chi(\overline M^2)+n
\end{equation}
 holds,
 where $K_{d\sigma^2}$ {\rm(}resp.~$dA_{d\sigma^2}${\rm)} are the Gaussian
 curvature {\rm (}resp.\ the area element{\rm)} of $d\sigma^2$, and 
 $\chi(\cdot)$ denotes the Euler number.
 The equality of \eqref{eq:osserman-curvature} holds if and only if 
 all ends are regular and cylindrical.
\end{proposition}
\begin{proof}
 Since $d\sigma^2$ is complete and $K_{d\sigma^2}$ is non-positive,
 \eqref{eq:osserman-curvature} is nothing but the Cohn-Vossen inequality
 for complete Riemannian $2$-manifold \cite{CV}.
 More precisely,
\[
   \frac{1}{2\pi}\int_{M^2}(-K_{d\sigma^2})\,dA_{d\sigma^2}
     =-\chi(\overline M^2)
       +\sum_{j=1}^n \ord_{p_j} d\sigma^2
\]
 holds (see \cite{Shiohama} or \cite{Fang}).
 Here, since $d\sigma^2$ is complete, $\ord_{p_j}d\sigma^2\leq -1$.
 Hence the equality in \eqref{eq:osserman-curvature} holds if and only
 if $\ord_{p_j}d\sigma^2=-1$ hold for all $j=1,\dots,n$.
 On the other hand, by \eqref{eq:sasaki-metric},
 it holds that
\[
    \ord_{p_j}d\sigma^2 = \min\{\ord_{p_j}|\omega|^2,\ord_{p_j}|\theta|^2\}.
\]
 Hence if $\ord_{p_j}d\sigma^2=-1$, the end $p_j$ is regular 
 because of \eqref{def_hopf} and Lemma~\ref{lemma:GMM}, and then 
 $p_j$ is cylindrical because of Lemma~\ref{lem:cylindrical}.
\end{proof}
Note that the left-hand side of \eqref{eq:osserman-curvature}
may not be an integer.

\section{Examples and a Classification}
\label{sec:example}
In this section, we investigate complete flat fronts of which all ends
are regular and embedded. 
We shall classify them when the number of ends is less than or 
equal to 3. 
We begin by reviewing known examples and their hyperbolic Gauss maps. 

\begin{example}[flat fronts of revolution]\label{ex:ffofrev}
 Let $\overline M^2$ denote the Riemann sphere $S^2=\C \cup \{ \infty \}$ 
 and consider a pair $(G, G_*)$ of meromorphic functions on 
 $\overline M^2$ defined by $G(z)=z$, $G_*(z)=\alpha z$ for some constant 
 $\alpha \in \R \setminus \{1 \}$. Define $M^2$ as follows: 
\begin{equation}
 M^2 := \begin{cases}
	 \overline M^2 \setminus \{ 0 \} & \text{ if } \alpha =0 \\
	 \overline M^2 \setminus \{ 0, \infty \} & \text{ otherwise. }  
	\end{cases}
\end{equation}
One can easily check that $M^2$ and $(G,G_*)$ satisfy the conditions 
\ref{item:rep-1} and \ref{item:rep-2} of Theorem~\ref{thm:rep}. 
Indeed, these data give a Legendrian immersion  
\begin{equation}
 E = \begin{pmatrix}
       z ^{\frac{-\alpha}{1- \alpha}} / c & 
      {c \alpha} z^{\frac{1}{1- \alpha}}/(1- \alpha) \\
      z ^{\frac{-1}{1- \alpha}} / c & 
      {c} z^{\frac{\alpha}{1- \alpha}} /(1- \alpha)
     \end{pmatrix} \text{ for some constant $c$}
\end{equation} 
and a resulting flat front $f:=EE^* \colon M^2 \to H^3$. 
\begin{figure}
\begin{center}
 \begin{tabular}{c@{\hspace{1.5cm}}c}
  \includegraphics[width=3cm]{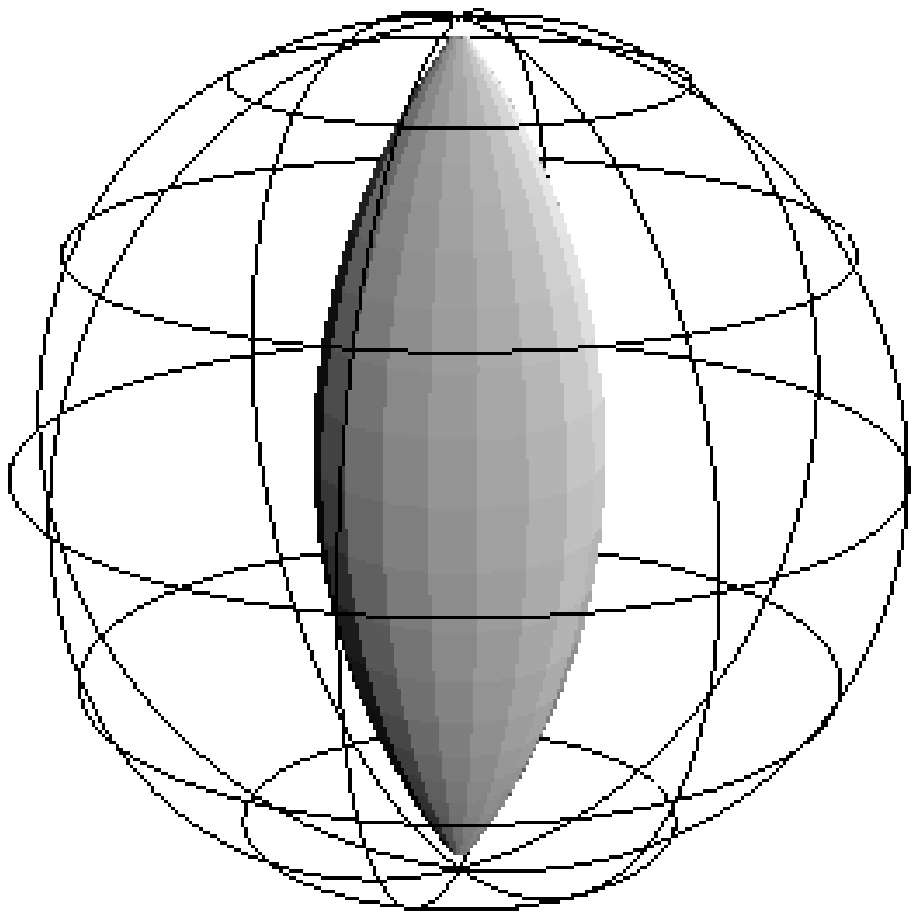} &
  \includegraphics[width=3cm]{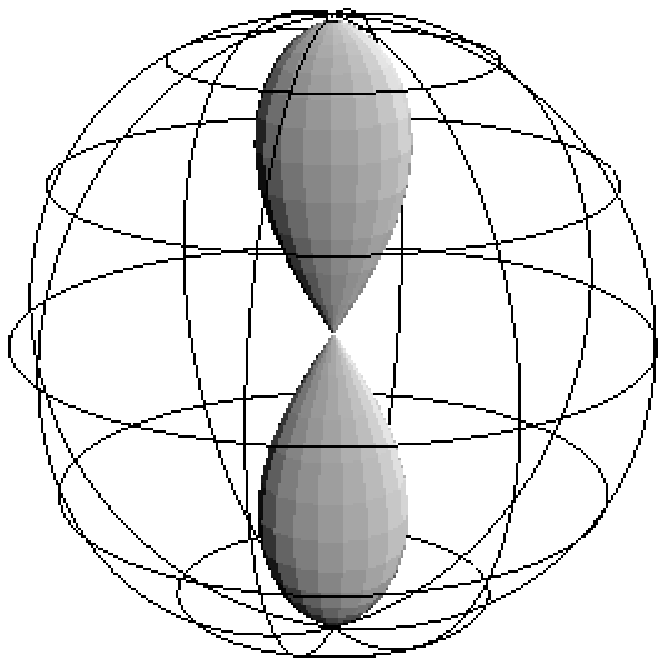} \\
   $\alpha=-1$ & $\alpha < 0$ ($\alpha\neq -1$) \\[12pt]
  \includegraphics[width=3cm]{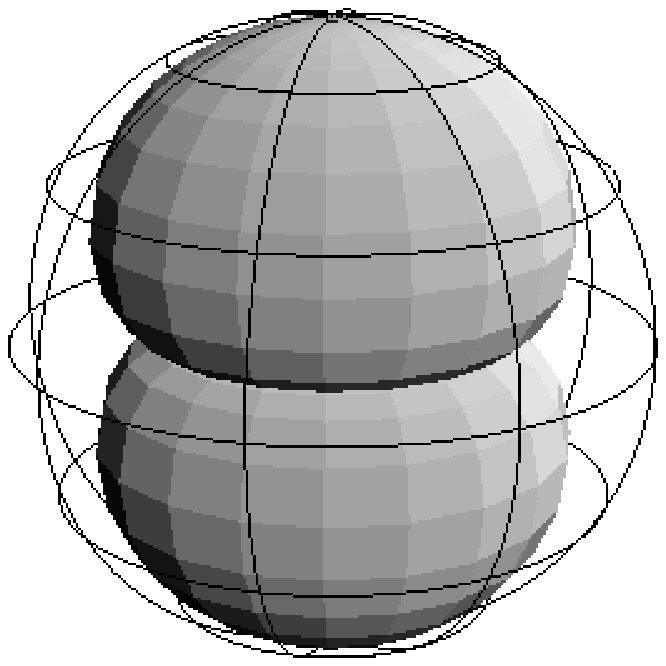} &
  \includegraphics[width=3cm]{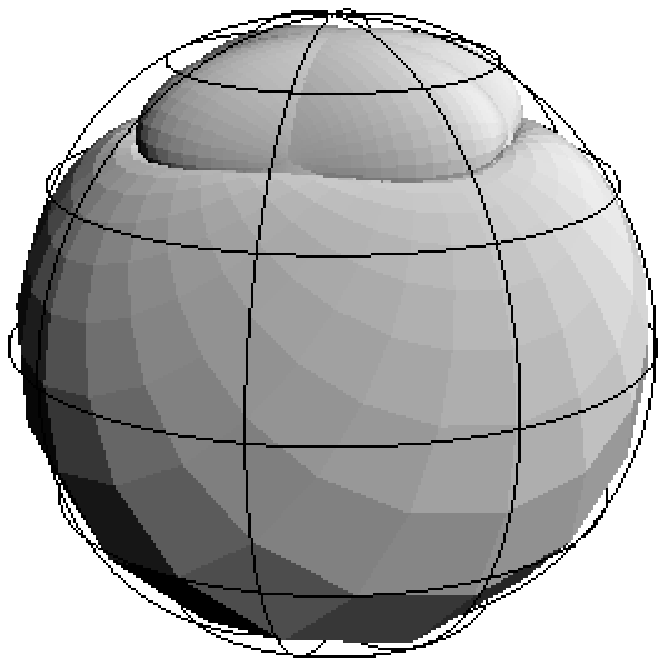} \\
  \multicolumn{2}{c}{$\alpha>0$}
 \end{tabular}
\end{center}
\caption{Flat fronts of revolution}
\end{figure}
This flat front $f$ is a horosphere if $\alpha =0$ or a hyperbolic
cylinder if $\alpha =-1$. 
We shall call $f$ an {\it hourglass\/}
if $\alpha (\ne -1) < 0$, or a {\it snowman\/} if $\alpha > 0$. 
The first and second canonical forms and the Hopf differential  are
represented as
\begin{gather*}
   \omega = -\frac{1}{c^2}z^{-2/(1-\alpha)}\,dz,\qquad
   \theta = \frac{c^2\alpha}{(1-\alpha)^2}z^{2\alpha/(1-\alpha)}\,dz\\
   Q = -\frac{\alpha}{(1-\alpha)^2}z^2\,dz^2.
\end{gather*}
The total curvature of the pull-back of the Sasakian metric
is calculated as
\[
   \frac{1}{2\pi}\int_{M^2}(-K_{d\sigma^2})\,dA_{d\sigma^2}
          = 2 \left|\frac{1+\alpha}{1-\alpha}\right|.
\]
\end{example}

Horospheres can be characterized by the hyperbolic Gauss maps.
\begin{proposition}\label{prop:horo-1}%
 Let $f\colon M^2 \to H^3$ be a complete flat front. 
 Assume that one of the hyperbolic Gauss maps $G$, $G_*$ is constant. 
 Then $f$ is a horosphere.
\end{proposition}
\begin{proof}
 It suffices to prove the case when $G_*$ is constant. 
 In this case, $G$ must be non-branched, 
 because $\mathcal G =(G,G_*)$ is an immersion.  
 On the other hand, $Q$ is identically zero.
 It follows from Lemma \ref{lemma:GMM} that $G$ has at most poles, that is, 
 $G$ is a meromorphic function on a compact 
 Riemann surface $\overline M^2$. 
 This implies that $G$ gives a biholomorphism  $\overline M^2 \cong S^2$. 
 Therefore we may assume $G(z) = z$ on 
 $\overline M^2 (\cong S^2 \cong \C \cup \{ \infty \})$. 
 Then it follows from Example~\ref{ex:ffofrev} that $f$ is a horosphere.   
\end{proof}
\begin{lemma}\label{lem:nlessthan3}
 Let $f\colon M^2=\overline M^2 \setminus \{p_1, \dots , p_n\} \to H^3$ 
 be a complete flat front with embedded regular ends $p_1, \dots p_n$. 
 If $n \le 3$, then $\overline M^2$ is biholomorphic to the Riemann sphere.   
\end{lemma}
\begin{proof}
 By Proposition \ref{prop:horo-1}, it suffices to prove this when  both
 $G$ and $G_*$ are non-constant, i.e.\ $\deg G \ge 1$ and 
 $\deg G_* \ge 1$. 
 Since all ends are regular and embedded, 
 $\deg G + \deg G_* = n \le 3$ holds. 
 Therefore $\deg G =1$ or $\deg G_* =1$. 
 This implies that $G$ or $G_*$ is a biholomorphism to the Riemann sphere.  
\end{proof}
Let us investigate complete flat fronts 
$f \colon \overline M^2 \setminus \{p_1, p_2\} \to H^3 $
with two embedded regular ends $p_1,p_2$. 
As stated in Lemma~\ref{lem:nlessthan3}, 
$\overline M^2 = S^2 \cong \C \cup \{ \infty \}$.  
Without loss of generality, we may assume that the images of two ends 
are $0, \infty \in S^2(= \partial H^3)$ respectively, that is, 
\begin{equation}\label{eq:imageofends}
 G(p_1)=G_*(p_1)=0 \text{ and } G(p_2)=G_*(p_2)= \infty.  
\end{equation}
It follows 
from the embeddedness of the ends that both $G$ and $G_*$ have 
degree $1$. 
We identify $\overline M^2$ with $S^2$ via $G$, that is, $G(z)=z$. 
Then the coordinates of $p_1$, $p_2$ are $z=0, \infty$, respectively. 
On the other hand, we can set $G_*(z)= (az+b)/(cz+d)$. 
It follows from \eqref{eq:imageofends} that $b=c=0$. 
Therefore $G_*(z) = \alpha z$ for some nonzero constant $\alpha$. 
Moreover, the conditions \ref{item:rep-1} and \ref{item:rep-2} of
Theorem~\ref{thm:rep} imply $\alpha \in \R \setminus \{0,1\}$. 

To summarize, $f$ is congruent to a flat front of  
$(G, G_*)=(z, \alpha z)$ for some $\alpha \in \R \setminus \{0,1\}$. 
Hence, it is a flat front of revolution (see Example~\ref{ex:ffofrev}).

\medskip

Next we investigate complete flat fronts 
$f \colon \overline M^2 \setminus \{p_1, p_2, p_3\} \to H^3 $, 
called {\it trinoids\/}, 
with $3$ embedded regular ends $p_1,p_2, p_3$. 
We may assume $\overline M^2 = S^2 \cong \C \cup \{ \infty \}$ by 
Lemma~\ref{lem:nlessthan3} and $\deg G=1$, $\deg G_* =2$. 
Similar to the case of two-end fronts above, we may assume that 
$G(z)=z$ and 
\begin{equation}\label{eq:imageof3ends}
 G_*(0) = 0, \quad G_*(1) = 1, \quad G_*(\infty) = \infty 
\end{equation} 
are the images of the ends. 
Since $G_*$ is a meromorphic function on $\overline M^2=S^2$ of degree
$2$, it is a fraction of polynomials of degree $ \le 2$. 
Indeed, it is verified from \eqref{eq:imageof3ends} 
that $G_*$ is the following form:
\begin{equation}\label{eq:Gstar-1}
 G_*(z)= \frac{z(\alpha z + \beta)}{ \gamma z+1},
\end{equation} 
where $\alpha,\beta,\gamma \in \C$ satisfy
\begin{equation}\label{eq:Gstar-2}
  \alpha+\beta = \gamma+1, \ \alpha \ne 0, \  
  \alpha-\beta\gamma \ne 0.  
\end{equation}
The conditions \eqref{eq:Gstar-2} can be rewritten as 
\begin{equation}\label{eq:Gstar-3}
  \beta = -\alpha+\gamma+1, \ \alpha(\alpha-\gamma)(\gamma+1) \ne 0.  
\end{equation} 

By straightforward computation, we can see that 
\begin{equation*}
 \frac{dG}{G-G_*}= \left(\frac{1}{\gamma-\alpha} \right)
 \frac{\gamma z+1}{z(z-1)}dz, 
\end{equation*} 
which has poles only at $z=0, 1, \infty$. 
All of them are simple poles of which residues are $-1/(\gamma-\alpha)$, 
$(\gamma+1)/(\gamma-\alpha)$, $-\gamma/(\gamma-\alpha)$, 
respectively. 
These residues must be real,
because of condition \ref{item:rep-2} of Theorem \ref{thm:rep}. 
Hence we have $\alpha, \gamma \in \R$ ($\beta \in \R$). 
It follows from Theorem \ref{thm:rep} and \eqref{eq:omega} 
that 
\begin{align*}
 \xi & \left(=c \exp \int \frac{dG}{G-G_*}\right)
 = c \, z^{\frac{1}{\alpha -\gamma}}
  (z-1)^{\frac{\gamma +1}{\gamma -\alpha}}, \\
 \omega & \left( = - \xi^{-2}dG \right)
 = -c^{-2}\, z^{\frac{2}{\gamma -\alpha}}
 (z-1)^{\frac{2\gamma+2}{\alpha-\gamma}}dz. 
\end{align*}
Furthermore, the Hopf differential $Q$ is computed as 
\begin{equation*}
 Q\left( = -\frac{dGdG_*}{(G-G_*)^2}\right) = 
 -\frac{1}{(\gamma-\alpha)^2}
\frac{\alpha\gamma z^2 +2\alpha z +\beta}{z^2(z-1)^2} dz^2. 
\end{equation*} 
Thus we see that $Q$ has poles only at $z=0,1,\infty$
with orders at most $2$.
Indeed, 
\begin{multline*}
 (\ord_0Q, \ord_1Q, \ord_{\infty}Q)\\
 =
\begin{cases}
 (-1,-2,-2) 
     & \text{if $\alpha=\gamma+1$ ($ \iff \beta=0$)} \\ 
 (-2,-1,-2) 
     & \text{if $\alpha=-1$ ($\iff \alpha\gamma+2\alpha+\beta=0$)} \\ 
 (-2,-2,-1) 
     & \text{if $\gamma=0$ ($\iff \alpha\gamma=0$)} \\ 
 (-2,-2,-2) & \text{otherwise.}
\end{cases}
\end{multline*} 
To summarize the arguments above, we obtain the following 
classification theorem. 

\begin{theorem}\label{thm:ff3end}
 Let $f\colon M^2 \to H^3$ be a complete flat front of which all ends 
 are regular and embedded. If the number of ends of $f$ is at most
 $3$, then $f$ is congruent to one of the following{\rm :}
\begin{enumerate}[\rm (i)]
 \item a horosphere if it has one end, 
 \item a hyperbolic cylinder, an hourglass, or a snowman 
       if it has two ends, 
 \item a trinoid with 
       $\displaystyle (G,G_*) =  
       \left(z,\frac{z(\alpha z + \beta)}{\gamma z+1} \right)$
       where $\alpha,\beta,\gamma$ are real constants satisfying 
       \eqref{eq:Gstar-3}, if it has three ends.
\end{enumerate}
\end{theorem}
In contrast to Theorem~\ref{thm:ff3end}, 
for arbitrary distinct points 
$p_1, \dots ,p_n \in \partial H^3 = \C \cup \left\{ \infty \right\}$,  
we can construct a flat front of genus zero with embedded regular ends 
$p_1, \dots ,p_n$ as follows: 
\begin{example}
 Let $p_1, \dots ,p_n$ be arbitrary distinct points in 
 $\partial H^3 = \C \cup \left\{ \infty \right\}$.   
 Without loss of generality, we may assume that $p_n=\infty$. 
 Let us choose non-zero real numbers $a_1, \dots, a_{n-1}$   
 so that $a_1 + \dots + a_{n-1} \ne 0, 1$. 
 We set
\begin{align*}
    M^2 &= \C\setminus\{p_1,\dots,p_{n-1}\}, \\
    G   &= z ,\\
    G_* &= \left. \left(
             z \sum_{k=1}^{n-1}
               \left\{a_k \prod_{j \ne k}\left(z-p_j\right)\right\}
             - \prod_{j=1}^{n-1}(z-p_j)\right) \right/ 
              {\sum_{k=1}^{n-1}\left\{a_k \prod_{j \ne k}
                \left(z-p_j\right)\right\}}, 
\end{align*}
then 
\begin{equation}\label{eq:tau_nnoid}
 \frac{dG}{G-G_*} =  
 \left\{ \frac{a_1}{z-p_1}+\frac{a_1}{z-p_1}+ \dots + 
 \frac{a_{n-1}}{z-p_{n-1}} \right\}dz.  
\end{equation}
It follows from \eqref{eq:tau_nnoid} that $M^2$ and $(G,G_*)$ satisfy
the assumptions of Theorem~\ref{thm:rep}. 
Therefore these data yield a flat front. 
By straightforward computation, it is easily seen that 
\begin{align*}
 \xi &= c \prod_{j=1}^{n-1}(z-p_j)^{a_j}, \\ 
 \omega & \left(= -dG/\xi^2\right)=
        -c^{-2} \left\{ \prod_{j=1}^{n-1}(z-p_j)^{-2a_j} \right\}dz, \\
 Q &= \left\{ \sum_{j=1}^{n-1} \frac{a_j}{(z-p_j)^2} 
    - \left( \sum_{j=1}^{n-1} \frac{a_j}{z-p_j}\right)^2\right\}dz^2. 
\end{align*}
It follows that $p=p_1, \dots, p_{n-1}, \infty$ are complete regular 
ends, respectively. 
Moreover, they are embedded ends since 
$\deg G + \deg G_*\left(=1+(n-1)=n\right)$ is equal to 
the number of ends (see Theorem~\ref{thm:Oss-ineq}).
\end{example}
Finally, we give examples of a complete flat front of genus $1$.
\begin{figure}
\begin{center}
 \begin{tabular}{c@{\hspace{1.5cm}}c}
  \includegraphics[width=3cm]{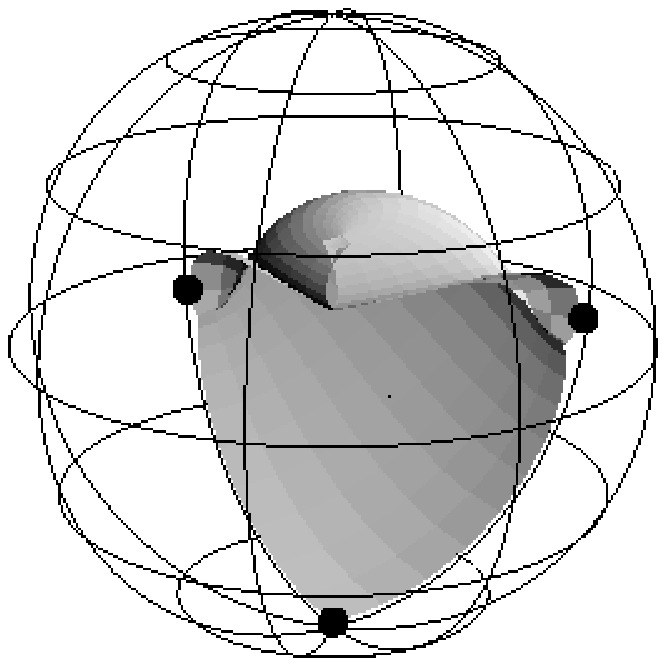} &
  \includegraphics[width=3cm]{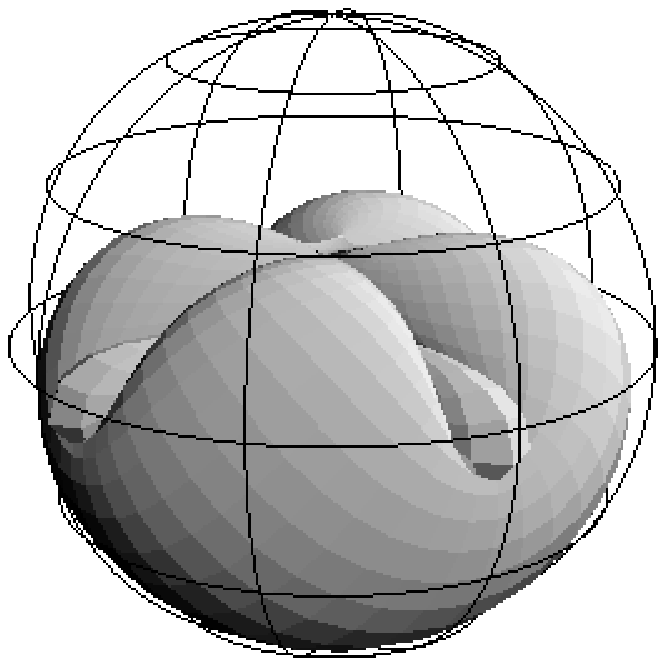} \\
   $c=0.3$ ($1/4$-cut) &  $c=0.3$ \\[12pt] 
  \includegraphics[width=3cm]{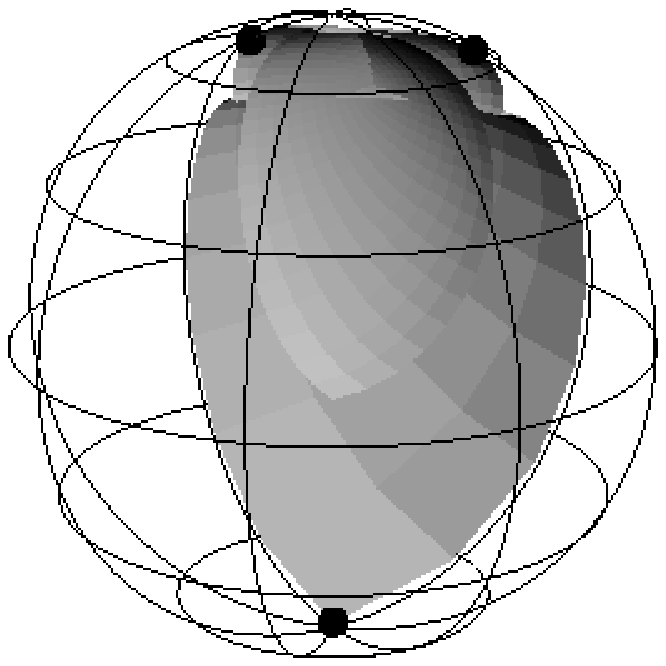} &
  \includegraphics[width=3cm]{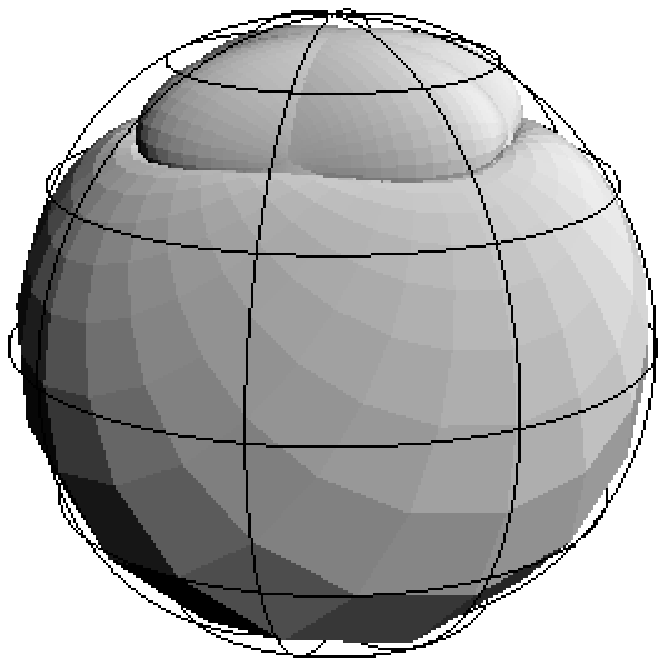} \\
   $c=1.0$ ($1/4$-cut) &  $c=1.0$
 \end{tabular}
\end{center}
\caption{%
 Genus one flat fronts with 5 embedded ends (Example~\ref{ex:torus}).
 In the figures of $1/4$-cut (left), the ends are shown as the dotted
 points.
}
\end{figure}
\begin{example}[of genus 1, with 5 embedded ends]
\label{ex:torus}
 Let $\wp$ denote the Weierstrass $\wp$ function on the 
 square torus $T^2=\C / \{\Z \oplus i\Z \}$. 
 We note that $\wp$ satisfies the following differential equation: 
\begin{equation*}
 \left( \wp' \right)^2 = 4 \wp \left(\wp^2 -e_1^2 \right), \qquad
 e_1= \wp(1/2). 
\end{equation*}
Take two meromorphic functions 
\begin{equation}\label{eq:torusGG}
 G=\wp', \quad G_* = - \frac{8e_1^2}{3} \frac{\wp}{\wp'}
\end{equation} 
on $T^2$. Let $M^2$ be a Riemann surface $T^2$ punctured at 
five points where $G$ and $G_*$ take the same value, i.e.  
\begin{equation}\label{eq:p_torus}
 M^2 := T^2 \setminus \{z \ ; \ \wp(3\wp^2-e_1^2)=0 \}. 
\end{equation}
We remark that $\wp$ has a double zero at $z=(1+i)/2$, 
and $3\wp^2-e_1^2$ has four simple zeros. 

For these data, a computation gives
\begin{equation*}
  \frac{dG}{G-G_*}=\frac{3}{2}\frac{\wp'}{\wp}. 
\end{equation*} 
This implies that the conditions \ref{item:rep-1} and
\ref{item:rep-2} of Theorem \ref{thm:rep} are satisfied. 
Therefore, the Riemann surface \eqref{eq:p_torus} 
and meromorphic functions \eqref{eq:torusGG} define a flat front. 

The first canonical form $\omega$ and the Hopf differential $Q$ are
computed to be 
\begin{align*}
 \omega &= -\frac{2}{c^2} \frac{3 \wp^2 -e_1^2}{\wp^3}dz, \\ 
 Q &= \frac{-6e_1^2(\wp^2+e_1^2)}{\wp(3\wp^2-e_1^2)} dz\,dz, 
\end{align*} 
from which the completeness of the ends 
$\{z \ ; \ \wp(3\wp^2-e_1^2)=0 \}$ follows. 
One can also verify the consistency of the data $G$, $\omega$ and $Q$ by
the formula (1.10) in \cite{KUY}.
Obviously, all ends are regular. 

On the other hand,  
since $G_*$ has only simple zeros at $z=0, (1+i)/2$, 
the degree of $G_*$ is equal to $2$, i.e., 
$\deg G_*=2$. 
Furthermore, it is obvious that $\deg G =3$. 
Hence, the equality in Theorem \ref{thm:Oss-ineq} is attained. 
Therefore all five ends are embedded. 
\end{example}

\end{document}